# A GAUSSIAN KINEMATIC FORMULA[1]

By Jonathan E. Taylor

*Stanford University*

In this paper we consider probabilistic analogues of some classical integral geometric formulae: Weyl–Steiner tube formulae and the Chern–Federer kinematic fundamental formula. The probabilistic building blocks are smooth, real-valued random fields built up from i.i.d. copies of centered, unit-variance smooth Gaussian fields on a manifold $M$. Specifically, we consider random fields of the form $f_p = F(y_1(p), \ldots, y_k(p))$ for $F \in C^2(\mathbb{R}^k; \mathbb{R})$ and $(y_1, \ldots, y_k)$ a vector of $C^2$ i.i.d. centered, unit-variance Gaussian fields.

The analogue of the Weyl–Steiner formula for such *Gaussian-related* fields involves a power series expansion for the Gaussian, rather than Lebesgue, volume of tubes: that is, power series expansions related to the marginal distribution of the field $f$. The formal expansions of the Gaussian volume of a tube are of independent geometric interest.

As in the classical Weyl–Steiner formulae, the coefficients in these expansions show up in a kinematic formula for the expected Euler characteristic, $\chi$, of the excursion sets $M \cap f^{-1}[u, +\infty) = M \cap y^{-1}(F^{-1}[u, +\infty))$ of the field $f$.

The motivation for studying the expected Euler characteristic comes from the well-known approximation $\mathbb{P}[\sup_{p \in M} f(p) \geq u] \simeq \mathbb{E}[\chi(f^{-1}[u, +\infty))]$.

**1. Introduction.** In this paper we consider perhaps the simplest non-Gaussian models of smooth random fields on a manifold $M$. Our fields, which we refer to as *Gaussian related*, are smooth, real-valued random fields built up from i.i.d. copies of centered, unit-variance smooth Gaussian fields on a manifold $M$. Given a $C^2$, centered, unit-variance Gaussian process $y$

Received November 2002; revised November 2004.

[1]Supported in part by Natural Sciences and Engineering Research Council and the Israel Science Foundation.

on $M$ and $F \in C^2(\mathbb{R}^k; \mathbb{R})$ we consider the field

$$f_p = F(y_1(p), \ldots, y_k(p)),$$

where the fields $y_i$ are i.i.d. copies of $y$.

For a concrete example: set

$$F(x) = \|x\|^2,$$

then the field

$$f(p) = F(y_1(p), \ldots, y_k(p))$$

has a $\chi_k^2$ marginal distribution, and has previously been referred to as a "$\chi_k^2$" field [1, 20]. Note that, unlike a Gaussian field, a field with $\chi_k^2$ marginal distributions is not determined solely by its covariance function. In this work the fields are characterized by their covariance functions and the function $F$.

The motivation for studying the expected Euler characteristic of the excursions $f^{-1}[u, +\infty) \subset M$ comes from the approximation

$$\mathbb{P}\left[\sup_{p \in M} f(p) \geq u\right] \simeq \mathbb{E}[\chi(f^{-1}[u, +\infty))]$$

[1, 2, 20, 22]. This approximation has found uses in medical imaging, astrophysics and multivariate analysis [15, 16, 20, 23]. A heuristic justification of the above approximation can be found in [2, 6], with a rigorous justification in [17].

Our main result, Theorem 4.1, expresses the expected Euler characteristic, $\chi$, of the excursion sets $f^{-1}[u, +\infty)$

$$\mathbb{E}[\chi(M \cap f^{-1}[u, +\infty))]$$

in terms of geometric quantities related to the curvature of the Riemannian structure induced by the covariance function of $y$ [18] and geometric quantities related to the curvature of $F^{-1}[u, +\infty) \subset \mathbb{R}^k$, viewed as a subset of the probability space $(\mathbb{R}^k, \gamma_{\mathbb{R}^k})$, where $\gamma_{\mathbb{R}^k}$ is the standard Gaussian measure on $\mathbb{R}^k$:

$$\gamma_{\mathbb{R}^k}(A) = \int_A (2\pi)^{-k/2} e^{-\|x\|^2/2} \, dx.$$

That is, the geometric quantities related to $F^{-1}[u, +\infty)$ depend on both the curvature of $F^{-1}[u, +\infty)$ as well as the measure $\gamma_{\mathbb{R}^k}$.

In Theorem 4.1, curvature enters as coefficients in certain volume of tubes expansions. As these expansions are slightly nonstandard, we recall some basic facts about such expansions, dating back to Weyl and Steiner [12, 14, 19]. These expansions give an expression for the volume of a tubular neighborhood of a set $M \subset \mathbb{R}^n$, assumed either to be convex or an embedded submanifold, in terms of certain intrinsic measures on $M$, the so-called Lipschitz–Killing curvature (signed) measures $(\mathcal{L}_j(M, \cdot))_{0 \leq j \leq n}$. These measures depend on $M$ and are finitely additive in $M$. When $M$ is an embedded manifold in $\mathbb{R}^n$, the tube formula is due to Weyl [19], and when $M$ is a compact, convex domain, the expansion is generally attributed to Steiner (cf. [14]). The two formulae were generalized to sets of positive reach by Federer [11] who defined the Lipschitz–Killing curvature measures of such a set.

Weyl and Steiner's formulae state that for $r$ small enough

$$\mathcal{H}_k(T(M, r)) \triangleq \mathcal{H}_k(\{y \in \mathbb{R}^k : d(y, M) \leq r\}) = \lambda_{\mathbb{R}^k}(\{y \in \mathbb{R}^k : d(y, M) \leq r\})$$

(1.1)

$$= \sum_{j=0}^{k} \mathcal{L}_j(M, \mathbb{R}^k) \omega_{k-j} r^{k-j} \triangleq \sum_{j=0}^{k} \mathcal{L}_j(M) \omega_{n-j} r^{k-j},$$

where $d(\cdot, M)$ is the standard distance function on $\mathbb{R}^n$; $\omega_k = \pi^{k/2} / \Gamma(k/2 + 1)$, the volume of the unit ball in $\mathbb{R}^k$; $\mathcal{H}_k$ is $k$-dimensional Hausdorff measure and $\lambda_{\mathbb{R}^k}$ is Lebesgue measure. It is sometimes useful to rewrite (1.1) in terms of *Minkowski functionals* defined by

$$\frac{\mathcal{M}_{k-j}(M)}{(k-j)!} = \mathcal{L}_j(M) \omega_{k-j}$$

so that (1.1) reads as a (finite) Taylor series expansion

$$\lambda_{\mathbb{R}^k}(T(M, r)) = \sum_{j=0}^{k} \mathcal{M}_j(M) \frac{r^j}{j!}.$$

One of the deep facts about the Lipschitz–Killing curvatures, first proven by Weyl, is that they are *intrinsic* to $M$. That is, if we embed $M$ into a different Euclidean space in an isometric fashion, the Lipschitz–Killing curvature measures of $S$ are unchanged. They are also *intrinsic* in the Riemannian sense, because they are local and can be computed from a Riemannian metric on $M$.

However, because $\mathcal{L}_j(M)$ are intrinsic to $M$, the Minkowski functionals of $M$ are not intrinsic. In particular they depend on the dimension of the Euclidean space in which $M$ is embedded. Strictly speaking, it is therefore necessary to write $\mathcal{M}_j^{\lambda_{\mathbb{R}^k}}(M)$ to clarify which Euclidean space we are talking about, where $\lambda_{\mathbb{R}^k}$ refers both to which dimension $M$ is considered to be embedded in as well as which measure we are using to compute the volume of the tube. In Section 3 we will compute the volume of a tube with measures other than Lebesgue; specifically we will compute a Taylor series expansion of the standard Gaussian volume $\gamma_{\mathbb{R}^k}$ of certain tubes. These expansions will play a key role in our analogue of the Chern–Federer kinematic fundamental formula (KFF) [7, 10, 11, 13] in which the Lipschitz–Killing curvatures also play a prominent role.

The KFF relates the "averaged" $j$th Lipschitz–Killing curvature $\mathcal{L}_j(M_1 \cap gM_2)$ to the Lipschitz–Killing curvatures of $M_1$ and $M_2$, averaged over "typical" rigid motions $g$. Specifically, for $M_1$ and $M_2$, two embedded submanifolds of $\mathbb{R}^k$, the following relation holds:

$$(1.2) \qquad \int_{G_k} \mathcal{L}_j(M_1 \cap gM_2) \, d\mu_n(g)$$
$$= \sum_{i=0}^{k-j} \begin{bmatrix} i+j \\ j \end{bmatrix} \begin{bmatrix} k \\ i \end{bmatrix}^{-1} \mathcal{L}_{j+i}(M_1)\mathcal{L}_{k-i}(M_2),$$

where

$$\begin{bmatrix} k \\ i \end{bmatrix} = \binom{k}{i} \frac{\omega_k}{\omega_i \omega_{k-i}},$$

and $G_k = \mathbb{R}^k \times O(k)$ is the group of rigid motions on $\mathbb{R}^k$ with Haar measure $\mu_k = \lambda_{\mathbb{R}^k} \times \tilde{\mu}_k$, the product of Lebesgue measure and the invariant probability measure on $O(k)$.

The Chern–Gauss–Bonnet theorem [9, 11] states that $\mathcal{L}_0(\cdot) = \chi(\cdot)$, the Euler characteristic, so that we can rewrite the case $j = 0$ in (1.2) as an "averaged," or "expected," Euler characteristic

$$(1.3) \qquad \int_{G_k} \chi(M_1 \cap gM_2) \, d\mu_k(g) = \sum_{i=0}^{k} \begin{bmatrix} k \\ i \end{bmatrix}^{-1} \mathcal{L}_i(M_1)\mathcal{L}_{k-i}(M_2),$$

keeping in mind that the measure $\mu_k$ is not a finite measure.

This brings us to the subject studied in our work here. In this work we study the expected Euler characteristic of certain random sets derived from smooth zero-mean, unit-variance Gaussian random fields on a $C^3$ manifold $M$. Our main results are Gaussian analogues of (1.1) and (1.3), given in Corollary 3.4 and Theorem 4.1.

Specifically, let $y = (y_1, \ldots, y_k) : M \to \mathbb{R}^k$ be a sequence of i.i.d. zero-mean unit-variance Gaussian random fields (satisfying certain regularity conditions specified in Theorem 4.1 below) on a $C^3$ manifold $M$. We derive, using techniques similar to [1, 3, 18, 21], an expression for $\mathbb{E}[\chi(M \cap y^{-1}D)]$, where $D$ is a suitable $C^2$ domain in $\mathbb{R}^k$. In particular, we show in Lemmas 2.5 and 2.6 that if $M$ is a $C^3$ $n$-manifold, with or without boundary, and for suitable $F \in C^2(\mathbb{R}^k; \mathbb{R})$

$$(1.4) \qquad \mathbb{E}[\chi(M \cap y^{-1}(F^{-1}[u, +\infty)))] = \sum_{j=0}^{n} \mathcal{L}_j(M)\tilde{\rho}_j(F, u),$$

where the functions $\widetilde{\rho}_j(F, u)$, for $j \geq 1$, are Gaussian analogues of $\mathcal{M}^{\lambda_{\mathbb{R}^k}}$ and are given by integrals of certain functions over $F^{-1}\{u\} = \partial(F^{-1}[u, +\infty)) \subset \mathbb{R}^k$. In (1.4), the quantities $(\mathcal{L}_j(M))_{0 \leq j \leq n}$ are the (total) Lipschitz–Killing curvature measures of the Riemannian manifold $(M, g)$, where $g$ is the Riemannian metric induced by the random field $f$

$$g(X_p, Y_p) = \mathbb{E}[X_p f Y_p f]$$

for tangent vectors $X_p$ and $Y_p$ in $T_p M$.

When the $y_i$'s are defined on $\mathbb{R}^n$ and are isotropic, the functionals $\widetilde{\rho}_j(F, u)$ agree with the EC densities $\rho_{f,j}(u)$ of $F \circ y$ (cf. [21]), where it is shown that for compact $C^2$ domains $S$

$$(1.5) \qquad \mathbb{E}[\chi(S \cap f^{-1}[u, +\infty))] = \sum_{j=0}^{n} a_{j,n} \mathcal{M}_{n-j}^{\lambda_{\mathbb{R}^n}}(S) \mu_2^{j/2} \rho_{f,j}(u),$$

where $\mu_2$ is a spectral parameter of $f$, namely the variance of its first partial derivatives and $a_{j,n}$ are constants, independent of $S$ or $f$. It should be noted that (1.5) holds for any smooth isotropic field on $\mathbb{R}^n$, not just Gaussian-related fields, that is, those derived from i.i.d. Gaussians. Lemmas 2.5 and 2.6 are thus extensions of (1.5) to fields derived from i.i.d. Gaussians on manifolds. Further, these lemmas show that in order to compute

$$\mathbb{E}[\chi(M \cap y^{-1}(F^{-1}[u, +\infty)))],$$

we need only compute $(\mathcal{L}_j(M))_{0 \leq j \leq n}$ and $\widetilde{\rho}_j(F, u)$.

Lemmas 2.5 and 2.6 also provide a direct way to compute $\widetilde{\rho}_j(F, u)$, and hence $\rho_{f,j}(u)$ for isotropic Gaussian-related fields of the form $F \circ y$, in that they are represented as conditional expectations of random variables defined on $\mathbb{R}^k$ and do not contain any information concerning the derivative of the fields themselves. This should be compared with earlier works (cf. [8, 20]) where $\rho_{f,j}(u)$ (assuming without loss of generality that $\mu_2 = 1$) is expressed as

$$(1.6) \qquad \rho_{f,j}(u) = \mathbb{E}[\det(\nabla^2 f_{|j,t}) | \nabla f_{|j,t} = 0, f_{|j,t} \geq u] \varphi_{\nabla f_{|j,t}}(0),$$

where $\nabla^2 f_{kl,t} = \partial^2 f(t)/\partial x_k \partial x_l$ is the Hessian of $f(t)$, $f_{|j}$ is the restriction of $f$ to a $j$-dimensional subspace of $\mathbb{R}^n$ and $\varphi_{\nabla f_{|j,t}}(0)$ is the density of $\nabla f_{|j,t}$ for some fixed $t$. In these works, it was therefore necessary to work out the full joint distribution of the field along with its first two derivatives at a given point $t$ and then carry out a fairly delicate conditioning argument.

The connection between the functionals $\widetilde{\rho}_j(F, u)$ and kinematic formulae is related to an extension of the Weyl–Steiner tube formulae. We discuss this extension in Section 3, where we derive a formal power series expression for the standard Gaussian volume of a tube around $D$, that is, we derive a formal expression for

$$(1.7) \qquad \int_{T(D,r)} (2\pi)^{-k/2} e^{-\|x\|^2/2} \, d\lambda_{\mathbb{R}^k}(x) = \sum_{j=0}^{\infty} \mathcal{M}_j^{\gamma_{\mathbb{R}^k}}(D) \frac{r^j}{j!}.$$

The notation $\mathcal{M}_j^{\gamma_{\mathbb{R}^k}}(D)$ is chosen to highlight the analogy between $\mathcal{M}_j^{\gamma_{\mathbb{R}^k}}(D)$, the coefficient of $r^j/j!$ in a formal expression for the Gaussian volume of $T(D,r)$ and $\mathcal{M}_j^{\lambda_{\mathbb{R}^k}}(D)$, the coefficient of $r^j/j!$ in a formal expression for the Lebesgue volume of $T(D,r)$. Both expansions are exact in some cases, for instance, when $D$ or $D^c$ is a compact $C^2$ domain.

The "take-home message" of this work is Theorem 4.1, where we link Sections 2 and 3. We prove, by direct computation, that for suitable $F \in C^2(\mathbb{R}^k; \mathbb{R})$,

$$\tilde{\rho}_j(F, u) = (2\pi)^{-j/2} \mathcal{M}_j^{\gamma_{\mathbb{R}^k}}(F^{-1}[u, +\infty)),$$

that is, the EC densities $\rho_{j,f}(u)$ for Gaussian-related isotropic fields are (up to constants) coefficients in a Gaussian volume of tubes expansion around $F^{-1}[u, +\infty)$. We can thus rewrite (1.4) as

$$(1.8) \qquad \mathbb{E}[\chi(M \cap y^{-1}D)] = \sum_{j=0}^{n} \mathcal{L}_j(M)(2\pi)^{-j/2} \mathcal{M}_j^{\gamma_{\mathbb{R}^k}}(D).$$

We conclude with some examples in Section 5, specifically rederiving, in light of (1.8), the EC densities of Gaussian and $\chi^2$ fields. While these are not new, their derivation sheds some light on the origin of the formulae. For instance, although it has been shown that the EC densities for a real-valued Gaussian field are essentially just Hermite polynomials times the standard Gaussian density, (1.8) shows that the reason the Hermite polynomials appear is the fact that they are derivatives of the Gaussian densities, hence (up to constants) the coefficients of powers of $r$ in a power series expansion of the Gaussian measure of the tube $[u-r, +\infty)$ of radius $r$ around $[u, +\infty)$.

The new results in Section 5 are the EC densities of the noncentral $\chi^2$ fields, as well as the EC densities of what we refer to as the "correlated conjunction" random field which is given by taking the minimum of two correlated Gaussian fields defined on a manifold $M$. Specifically given $\rho \in (-1, 1)$ and $y = (y_1, y_2)$ two independent centered, unit-variance Gaussian fields on $M$, define the new random fields

$$z_1 = y_1,$$

$$z_2 = \rho \cdot y_1 + \sqrt{1 - \rho^2} \cdot y_2,$$

$$z_1 \wedge z_2(p) = \min(z_1(p), z_2(p)).$$

In Section 5.4 we derive the EC densities of $z_1 \wedge z_2$ by relating the EC densities to the Gaussian measure of tubes around arbitrary cones in $\mathbb{R}^2$. The "correlated conjunctions" is a simple model for correlated random sets: the independent model ($\rho = 0$) has been used in certain neuroimaging applications [24].

The outline of the paper is as follows: Section 2 is dedicated to deriving an explicit expression for the EC densities $\rho_{f,j}(u) = \widetilde{\rho}_j(F, u)$. Section 3 is dedicated to extending the Weyl–Steiner tube formulae to measures with smooth densities with respect to Lebesgue measure and is possibly of independent interest. In Section 4 we show that the functionals $\widetilde{\rho}_j(F, u)$ are, up to constants, coefficients in a Taylor series expansion for the standard Gaussian measure of a tube around $F^{-1}[u, +\infty)$. Section 5 is devoted to examples.

## 2. Euler characteristic densities.
In this section we derive an explicit integral representation for the EC densities $\widetilde{\rho}_j(F, u)$. The proof involves some preliminary lemmas necessary to carry out the conditioning in (1.6).

### 2.1. *Gaussian random fields on manifolds.*
In this section we derive certain expectations for Gaussian random fields on a manifold $M$. We first review some linear algebraic preliminaries; readers are referred to [3, 18] for further details. For an $n$-dimensional vector space $V$, let $(\mathcal{T}^*(V), \otimes)$ denote the algebra of covariant tensors on $V$; $(\bigwedge^*(V), \wedge)$, the Grassman algebra and $(\bigwedge^*(V) \otimes \bigwedge^*(V), \cdot)$, the algebra of double forms $\bigwedge^*(V) \otimes \bigwedge^*(V) = \bigoplus_{r,s=0}^{\infty} \Lambda^{r,s}(V) = \bigoplus_{r,s=0}^{\infty} \Lambda^r(V) \otimes \Lambda^s(V)$ endowed with the "double-wedge" product

$$(\alpha \otimes \beta) \cdot (\gamma \otimes \delta) = (\alpha \wedge \gamma) \otimes (\beta \wedge \delta).$$

The subalgebra $(\bigoplus_{j=0}^{\infty} \Lambda^{j,j}(V), \cdot)$ is denoted by $(\bigwedge^{*,*}(V), \cdot)$. Any inner product $\langle \cdot, \cdot \rangle$ on $V$ induces an inner product on $\bigwedge^*(V)$, and the trace, $\mathrm{Tr} \colon \bigwedge^*(V) \otimes \bigwedge^*(V) \to \mathbb{R}$ is defined by

$$\mathrm{Tr}(\alpha \otimes \beta) = \langle \alpha, \beta \rangle,$$

and extended linearly.

We call $W \colon (\Omega, \mathcal{F}, \mathbb{P}) \to \Lambda^{1,1}(V)$ a *Gaussian double form* if, for any basis $B_V = \{v_1, \ldots, v_n\}$ of $V$, the matrix $W(v_i, v_j)$ has (jointly) Gaussian entries.

We recall the following useful lemma from [18].

LEMMA 2.1. *Suppose that $W$ is a Gaussian double form; then*

$$\mathbb{E}[W^k] = \sum_{j=0}^{\lfloor k/2 \rfloor} \frac{k!}{(k-2j)!\,j!\,2^j} \mu^{k-2j} C^j,$$

*where*

$$\mu \triangleq \mathbb{E}[W] \in \Lambda^{1,1}(V), \qquad C \triangleq \mathbb{E}[(W - \mathbb{E}[W])^2] \in \Lambda^{2,2}(V).$$

The next result we will need in subsequent sections concerns the conditional expectation of certain random double forms on finite-dimensional Hilbert spaces. Specifically, let $(V, \langle \cdot, \cdot \rangle)$ be a finite-dimensional Hilbert space and $L \subset V$, a subspace of $V$ with $P_L$ denoting orthogonal projection onto $L$. Suppose now that $(X_i)_{1 \leq i \leq N \leq \dim(V)}$ are i.i.d. $V$-valued random variables with common distribution $\gamma_V$, the canonical Gaussian measure on $V$

$$\gamma_V(A) = (2\pi)^{-\dim(V)/2} \int_A e^{-\|v\|^2/2} \, d\mathcal{H}_{\dim(V)}(v).$$

We evaluate the expression

(2.1)     $\mathbb{E}[\alpha((X_{i_1}, \ldots, X_{i_k}), (X_{j_1}, \ldots, X_{j_k})) | P_L X_1, \ldots, P_L X_N],$

for $\alpha \in \Lambda^{k,k}(V)$ and $k$-tuples $\tilde{i}_k = (i_1, \ldots, i_k)$ and $\tilde{j}_k = (j_1, \ldots, j_k)$ of $\{1, \ldots, N\}$.

Before evaluating (2.1), we recall the definition of the annihilation operator on $\mathcal{T}^*(V)$, as well that of a contraction on $\Lambda^*(V) \otimes \Lambda^*(V)$. For $v \in V$, the annihilation operator $i_v$ on $\bigwedge^*(V)$ by defining it on $\Lambda^r(V)$, as follows

$$(i_v \beta)(v_1, \ldots, v_{r-1}) = \beta(v, v_1, \ldots, v_{r-1})$$

and, if $\beta \in \Lambda^0(V)$, we set $i_v \beta = 0$. Although it is referred to as an annihilation operator, its real effect is to fix the value of the first coordinate to be $v$ when the form $i_v \beta$ is evaluated on $r-1$ other vectors. Since $\Lambda^*(V) \otimes \Lambda^*(V)$ consists of two "copies" of $\Lambda^*(V)$ the annihilation operator can act on either the first or second copy. To distinguish between these two cases we define operators $\eta_v$ and $\eta'_v$ on $\Lambda^*(V) \otimes \Lambda^*(V)$ by setting

$$\eta_v(\beta \otimes \gamma) = i_v \beta \otimes \gamma,$$
$$\eta'_v(\beta \otimes \gamma) = \beta \otimes i_v \gamma,$$

and extending linearly.

For $\alpha \in \bigwedge^*(V) \otimes \bigwedge^*(V)$

(2.2)
$$\eta_{v_1} \eta_{v_2} \alpha = -\eta_{v_2} \eta_{v_1} \alpha,$$
$$\eta'_{v_1} \eta'_{v_2} \alpha = -\eta'_{v_2} \eta'_{v_1} \alpha,$$
$$\eta'_{v_1} \eta_{v_2} \alpha = \eta_{v_2} \eta'_{v_1} \alpha,$$

as is easy to check from the definition. Also

$$\alpha((X_{i_1}, \ldots, X_{i_k}), (X_{j_1}, \ldots, X_{j_k})) = \left( \prod_{l=1}^k \eta_{X_{i_l}} \eta'_{X_{j_l}} \right) \alpha$$

so that the evaluation of (2.1) is really the evaluation of the conditional expectation of the polynomial

$$\prod_{l=1}^{k} \eta_{X_{i_l}} \eta'_{X_{j_l}} \in L(\Lambda^*(V), \Lambda^*(V)),$$

a random element of the set of linear maps from $\Lambda^*(V) \to \Lambda^*(V)$.

Fix $\{v_1, \ldots, v_{\dim(V)}\}$, an orthonormal basis for $(V, \langle \cdot, \cdot \rangle)$ and $L \subset V$ a subspace with orthonormal basis $\{v'_1, \ldots, v'_l\}$. We define the contraction operators $C, C_L, C_L^\perp$ on $\Lambda^*(V) \otimes \Lambda^*(V)$ as follows:

$$C\alpha = \sum_{i=1}^{\dim(V)} \eta_{v_i} \eta'_{v_i} \alpha, \qquad C_L \alpha = \sum_{i=1}^{l} \eta_{v'_i} \eta'_{v'_i} \alpha,$$

$$C_L^\perp \alpha = (C - C_L)\alpha = C_{L^\perp} \alpha,$$

where $L^\perp$ is the orthogonal complement of $L$ in $V$. For any two subspaces $L_1, L_2$, (2.2) implies

$$C_{L_1} C_{L_2} = C_{L_2} C_{L_1}$$

and for $\alpha \in \Lambda^{k,k}(V)$

$$C_L^k(\alpha) = k! \operatorname{Tr}^L(\alpha_{|L}),$$

where $\alpha_{|L}$ is the restriction of $\alpha$ to $L$.

With the notation established, we proceed to evaluate the conditional expectation (2.1).

LEMMA 2.2. *Suppose* $\alpha \in \bigwedge^*(V) \otimes \bigwedge^*(V)$, $X_i$, $1 \le i \le N$, *are i.i.d.* $V$-*valued random variables with distribution* $\gamma_V$, *and* $L_1, \ldots, L_N$ *are subspaces of* $V$. *Suppose further that* $\tilde{i}_p = (i_1, \ldots, i_p)$ *and* $\tilde{j}_q = (j_1, \ldots, j_q)$ *are two* $p$- *and* $q$-*tuples of distinct elements of* $\{1, \ldots, N\}$, *arranged such that the first* $n(\tilde{i}_p, \tilde{j}_q) \le \min(p, q)$ *elements of* $\tilde{i}_p$ *match those of* $\tilde{j}_q$ *and there are no further matches in the remaining elements of* $\tilde{i}_p$ *and* $\tilde{j}_q$. *The following relation then holds:*

$$\mathbb{E}\left[\left(\prod_{l=1}^{p} \eta_{X_{i_l}}\right)\left(\prod_{l=1}^{q} \eta'_{X_{j_l}}\right)\alpha \,\Big|\, P_{L_1} X_1, \ldots, P_{L_N} X_N\right]$$

$$(2.3) \qquad = \left(\prod_{l=1}^{n(\tilde{i}_p, \tilde{j}_q)} (C_{L_{i_l}}^\perp + \eta_{P_{L_{i_l}} X_{i_l}} \eta'_{P_{L_{i_l}} X_{i_l}})\right)\left(\prod_{m=n(\tilde{i}_p, \tilde{j}_q)+1}^{p} \eta_{P_{L_{i_m}} X_{i_m}}\right)$$

$$\times \left(\prod_{m=n(\tilde{i}_p, \tilde{j}_q)+1}^{q} \eta'_{P_{L_{i_m}} X_{i_m}}\right)\alpha.$$

REMARK. The restriction that the indices in $\tilde{\tilde{i}}_p$ and $\tilde{\tilde{j}}_q$ be distinct is to rule out the trivial case when either two of the elements of $\tilde{\tilde{i}}_p$ or $\tilde{\tilde{j}}_q$ are equal and the polynomial

$$\left(\prod_{l=1}^{p} \eta_{X_{i_l}}\right)\left(\prod_{l=1}^{q} \eta'_{X_{j_l}}\right) \equiv 0.$$

PROOF OF LEMMA 2.2. We consider first the case in which $p = q = 1$, and note that, by linearity, it suffices to prove the lemma for $\alpha \in \Lambda^r(V) \otimes \Lambda^s(V)$. Working from the definition of $\eta_{X_{i_1}} \eta'_{X_{j_1}}$

$$\eta_{X_{i_1}} \eta'_{X_{j_1}}(\alpha)((w_1, \ldots, w_{r-1}), (w'_1, \ldots, w'_{s-1}))$$

$$= \sum_{k,l=1}^{n} \langle X_{i_1}, v_k \rangle \langle X_{j_1}, v_l \rangle \alpha((v_k, w_1, \ldots, w_{r-1}), (v_l, w'_1, \ldots, w'_{s-1})),$$

where $B = \{v_1, \ldots, v_{\dim(V)}\}$ is some orthonormal basis for $(V, \langle \cdot, \cdot \rangle)$. If $i_1 = j_1$, we can further choose the orthonormal basis for $V$ so that $\{v_1, \ldots, v_{\dim(L_{i_1})}\}$ is an orthonormal basis for $L_{i_1}$ and $\{v_{\dim(L_{i_1})+1}, \ldots, v_{\dim(V)}\}$ is an orthonormal basis for $L_{i_1}^{\perp}$. The result then follows from well-known results about the conditional distributions of Gaussian vectors.

Next, consider $p, q \geq 1$,

$$\mathbb{E}\left[\left(\prod_{l=1}^{p} \eta_{X_{i_l}}\right)\left(\prod_{l=1}^{q} \eta'_{X_{j_l}}\right)\alpha \,\Big|\, P_{L_1} X_1, \ldots, P_{L_N} X_N\right]$$

$$= \mathbb{E}\left[\mathbb{E}\left[\left(\prod_{l=1}^{p} \eta_{X_{i_l}}\right)\left(\prod_{l=1}^{q} \eta'_{X_{j_l}}\right)\alpha \,\Big|\, X_1, \ldots, P_{L_{i_1}} X_{i_1}, \ldots, P_{L_{j_1}} X_{j_1}, \ldots, X_N\right]\Big|$$

$$P_{L_1} X_1, \ldots, P_{L_N} X_N\right]$$

$$= \mathbb{E}\left[\left(\prod_{l=2}^{p} \eta_{X_{i_l}}\right)\left(\prod_{l=2}^{q} \eta'_{X_{j_l}}\right)\right.$$

$$\left. \times \mathbb{E}[\eta_{X_{i_1}} \eta'_{X_{j_1}} \alpha | X_1, \ldots, P_{L_{i_1}} X_{i_1}, \ldots, P_{L_{j_1}} X_{j_1}, \ldots, X_N]\Big|\right.$$

$$P_{L_1} X_1, \ldots, P_{L_N} X_N\right]$$

$$= (\delta_{i_1 j_1} C_{L_{i_1}}^{\perp} + \eta_{P_{L_{i_1}} X_{i_1}} \eta'_{P_{L_{j_1}} X_{j_1}})$$

$$\times \mathbb{E}\left[\left(\prod_{l=2}^{p} \eta_{X_{i_l}}\right)\left(\prod_{l=2}^{q} \eta'_{X_{j_l}}\right)\alpha \,\Big|\, P_{L_1} X_1, \ldots, P_{L_k} X_k\right];$$

the general case then follows by iteration.  □

COROLLARY 2.3.  *Suppose $\alpha \in \Lambda^{k,k}(V)$, $X_1, \ldots, X_k$ are i.i.d. with distribution $\gamma_V$, $\tilde{i}_k$ and $\tilde{j}_k$ are as in Lemma* 2.2 *and $v_0$ is a unit vector in $V$. Then*

$$\mathbb{E}[\alpha((X_{i_1}, \ldots, X_{i_k}), (X_{j_1}, \ldots, X_{j_k}))|\langle X_1, v_0 \rangle, \ldots, \langle X_k, v_0 \rangle]$$

$$= \mathbb{1}_{\{n(\tilde{i}_k, \tilde{j}_k)=k\}} \left( k! \operatorname{Tr}^{v_0^\perp}(\alpha_{|v_0^\perp}) + \sum_{l=1}^{k} \eta_{\langle X_{i_l}, v_0 \rangle v_0} \eta'_{\langle X_{i_l}, v_0 \rangle v_0} C_{v_0^\perp}^{k-1}(\alpha) \right)$$

$$+ \mathbb{1}_{\{n(\tilde{i}_k, \tilde{j}_k)=k-1\}} (\eta_{\langle X_{i_k}, v_0 \rangle v_0} \eta'_{\langle X_{j_k}, v_0 \rangle v_0} C_{v_0^\perp}^{k-1}(\alpha)).$$

PROOF.  This is just an application of Lemma 2.2 which uses the fact that, if the two sets $\tilde{i}_k$ and $\tilde{j}_k$ differ by more than two indices, then the polynomial resulting from the application of Lemma 2.2 will be identically zero, since the term $\eta_v$ (as well as $\eta'_v$) will appear at least twice [cf. (2.2)].  □

For our purposes, we can reformulate the above corollary as follows.

COROLLARY 2.4.  *Suppose $y = (y_1, \ldots, y_j) : M \to \mathbb{R}^j$ are i.i.d. mean-zero, unit-variance Gaussian random fields on a smooth manifold $M$, $\nu$ is a vector field on $\mathbb{R}^j$ and $\alpha \in \Lambda^{k,k}(\mathbb{R}^j)$. Then, for any (nonrandom) set of vector fields $(Z_i)_{1 \le i \le m}$ on $M$, where $k \le m \le \dim(M)$*

$$\mathbb{E}[y^* \alpha | y, \langle y_* Z_i, \nu \rangle, 1 \le i \le m]$$

*is a random double form, whose restriction to $L = \operatorname{span}\{Z_1, \ldots, Z_m\}$ satisfies, for each $p \in M$,*

$$\mathbb{E}[y^* \alpha_{|L} | y, \langle y_* Z_i, \nu \rangle, 1 \le i \le m]$$

$$= \operatorname{Tr}^{\nu^\perp(y)}(\alpha(y)_{|\nu^\perp(y)}) I_{|L}^k + O\left( \sum_{i=1}^{m} \frac{(\langle y_* Z_i, \nu \rangle)^2}{\|\nu(y)\|^2} \|\alpha(y)\|_{\otimes^{2k} T_y \mathbb{R}^j} \right),$$

*in the sense that for each $p \in M$, on the set $\{\|\nu(y_p)\| > 0\}$*

$$|\mathbb{E}[y^* \alpha_{|L}((Z_{i_1}, \ldots, Z_{i_k}), (Z_{j_1}, \ldots, Z_{j_k}))_p | y_p, \langle y_* Z_i, \nu \rangle_p, 1 \le i \le m]$$

$$- \operatorname{Tr}^{\nu^\perp(y_p)}(\alpha(y_p)_{|\nu^\perp(y_p)}) I_{|L}^k((Z_{i_1}, \ldots, Z_{i_k}), (Z_{j_1}, \ldots, Z_{j_k}))_p|$$

$$\le K_{k,j} \left( \sum_{i=1}^{m} \frac{(\langle y_* Z_i, \nu \rangle_p)^2}{\|\nu(y_p)\|^2} \|\alpha(y_p)\|_{\otimes^{2k} T_{y_p} \mathbb{R}^j} \right),$$

*for some universal constant $K_{k,j}$ and all $k$-tuples $(i_1, \ldots, i_k)$ and $(j_1, \ldots, j_k)$.*

PROOF.    Without loss of generality, we can assume that the $(Z_i)_{1 \le i \le m}$ are orthonormal with respect to the metric induced by any one of the $y_i$'s. The result then follows from the fact that, for $p$ fixed, $(y_* Z_i(p))_{1 \le i \le n}$ are, conditional on $y(p)$, i.i.d. Gaussian random vectors in $T_{y(p)} \mathbb{R}^j$ with distribution $\gamma_{T_{y(p)} \mathbb{R}^j}$, combined with the previous corollary, along with the fact that for any vector field $W$ on $\mathbb{R}^j$,

$$|C^{k-1}_{\nu^\perp(y)} \eta_{P_{\eta(y)} W(y)} \eta'_{P_{\eta(y)} W(y)} \alpha(y)| < K'_{k,j} \frac{\langle W(y), \eta(y) \rangle^2}{\|\eta(y)\|^2} \|\alpha(y)\|_{\otimes^{2k} \mathbb{R}^j},$$

for some universal constant $K'_{k,j}$.    □

2.2.  *Integral representation of $\widetilde{\rho}_j(F, u)$.*    In this section we derive expressions for the expected Euler characteristic of the excursions of *Gaussian-related fields*, which we defined to be random fields of the form

$$f(p) = F \circ y(p),$$

where $y = (y_1, \ldots, y_k) : M \to \mathbb{R}^k$ are i.i.d. zero-mean unit-variance Gaussian random fields and $F \in C^2(\mathbb{R}^k)$. Recall that the EC densities were defined [21] for isotropic fields $f$ on $\mathbb{R}^n$ as follows:

$$(2.4) \qquad \mathbb{E}[\chi(M \cap f^{-1}[u, +\infty))] = \sum_{j=0}^{n} a_{j,n} \mathcal{M}_{n-j}^{\lambda_{\mathbb{R}^n}}(M) \mu_2^{j/2} \rho_{f,j}(u)$$

where the spectral parameter $\mu_2$ is just the variance of the first-order partial derivatives of $f$ and $\rho_{f,j}(u)$ depend only on the finite-dimensional distributions of $f$ and its first two derivatives at any spatial location. Suppose $f_1 = F \circ y_1$ and $f_2 = F \circ y_2$ are two (not identically distributed) isotropic Gaussian-related fields with $\mu_{2,1} = \mu_{2,2}$, that is, with equal variance of the first derivatives. Then, by inspecting the results in [20], it can be shown that the expected Euler characteristics, $f_i^{-1}[u, +\infty)$, will be identical (this will in fact be proven in the next lemma). In particular, the EC densities are functionals of $F$, so that, for isotropic Gaussian-related fields we can rewrite (2.4) as

$$(2.5) \qquad \mathbb{E}[\chi(M \cap f^{-1}[u, +\infty))] = \sum_{j=0}^{n} a_{j,n} \mathcal{M}_{n-j}(M) \mu_2^{j/2} \widetilde{\rho}_j(F, u)$$

for some functionals $\widetilde{\rho}_j$, $0 \le j < \infty$.

The simple case $k = 1$ and $F(x) = x$ or $F = \mathrm{Id}$, the identity map, leads to a real-valued Gaussian field, for which it is known [1] that

$$\widetilde{\rho}_j(\mathrm{Id}, u) = (2\pi)^{-(j+1)/2} \int_u^\infty H_j(x) e^{-x^2/2} \, dx,$$

with $H_j$ the $j$th Hermite polynomial. In the manifold setting, if $f$ is a real-valued zero-mean unit-variance Gaussian field on a $C^3$ manifold $M$, it was shown in [18] that

$$\mathbb{E}[\chi(M \cap f^{-1}[u, +\infty))] = \sum_{j=0}^{n} \mathcal{L}_j(M)\widetilde{\rho}_j(\mathrm{Id}, u),$$

that is, that (2.5) holds with $a_{j,n}\mathcal{M}_{n-j}(M)\lambda^{j/2}$ replaced by $\mathcal{L}_j(M)$, where the Riemannian metric with respect to which the $\mathcal{L}_j(M)$ are computed is the one induced by $f$, that is,

$$g_p(X_p, Y_p) = \mathbb{E}[X_p f Y_p f]$$

which is assumed to be $C^2$.

One way to think of the relation (2.5) is as a partitioning of the geometric information from the parameter space and information depending on the distribution of the field. The result in the real-valued Gaussian case suggests that this partitioning might hold for Gaussian-related fields on manifolds. This would mean that in order to compute expected Euler characteristics, we would simply need to calculate the Lipschitz–Killing curvatures of $(M, g)$ (recall that $g$ is the metric induced by any one of the coordinate random fields $y_i$) and the EC densities of a corresponding Gaussian-related isotropic field. The following lemma shows that this is indeed the case, and our work is reduced to computing the EC densities of isotropic Gaussian-related fields.

We set $\nabla f_E(p) = (E_1 f(p), \ldots, E_n f(p))$ for some $C^1$ section of $\mathcal{O}(M)$, the bundle of orthonormal frames on $(M, g)$ and

$$\alpha_\varepsilon(y) = (\omega_n \varepsilon^n)^{-1} \mathbb{1}_{B_{\mathbb{R}^n}(0, \varepsilon)}(y) \left( \bigwedge_{j=1}^{n} dy_j \right),$$

an approximation to the Dirac delta, as in [18].

LEMMA 2.5. Given $y = (y_1, \ldots, y_k)$ i.i.d. suitably regular centered, unit-variance Gaussian random fields on $M$, a $C^3$ $n$-manifold with or without boundary, and $F \in C^2(\mathbb{R}^k)$ is such that $f = F \circ y$ satisfies:

(i) $\lim_{\varepsilon \to 0} \mathbb{E}[\int_M \nabla f_E^*(\alpha_\varepsilon) \mathbb{1}_{\{f \geq u\}}] = \int_M \lim_{\varepsilon \to 0} \mathbb{E}[\nabla f_E^*(\alpha_\varepsilon) \mathbb{1}_{\{f \geq u\}}]$,

(ii) every critical point of $f$ in $f^{-1}(u - \delta, u + \delta)$ for some $\delta > 0$ is non-degenerate (cf. [18]), $\mathbb{P}$-a.s.

Then

$$\mathbb{E}[\chi(M \cap f^{-1}[u, +\infty))] = \sum_{j=0}^{n} \mathcal{L}_j(M)\widetilde{\rho}_j(F, u)$$

for some functionals $\widetilde{\rho}_j(F, u)$.

PROOF. We will prove the statement when $M$ has no boundary, omitting the calculations for $M$ with boundary, as they are similar to those in Theorem 5.1 of [18]. Further details can be found in [3].

Following the arguments in Theorem 4.1 of [18], under the assumptions of the theorem,

$$\mathbb{E}[\chi(M \cap f^{-1}[u, +\infty))]$$
$$= \int_M \lim_{\varepsilon \to 0} \mathbb{E}[\nabla f_E^*(\alpha_\varepsilon) \mathbb{1}_{\{f \geq u\}}]$$
$$= \frac{1}{n!} \int_M \lim_{\varepsilon \to 0} \mathbb{E}[\mathbb{1}_{\{f \geq u, \|\nabla f_E\| < \varepsilon\}} \operatorname{Tr}^M((-\nabla^2 f)^n)] \operatorname{Vol}_{M,g}$$
$$= \frac{1}{n!} \int_M \lim_{\varepsilon \to 0} \mathbb{E}[\mathbb{1}_{\{f \geq u, \|\nabla f_E\| < \varepsilon\}} \mathbb{E}[\operatorname{Tr}^M((-\nabla^2 f)^n) | f, \nabla f_E]] \operatorname{Vol}_{M,g}.$$

We will therefore devote our attention to evaluating

$$\mathbb{E}[(-\nabla^2 f)^n | f, \nabla f_E] \in \Lambda^{n,n}(M),$$

which we will show, for each $p \in M$ is of the form

$$(2.6) \qquad \mathbb{E}[(-\nabla^2 f)^n | f, \nabla f_E]_p = \sum_{j=0}^{\lfloor n/2 \rfloor} \frac{1}{j!} R_p^j A_j(p)$$

for some random double forms $A_j(p)$ [measurable with respect to $\sigma(f(p), \nabla f_E(p))$]. Since the distribution of the random vector $(f(p), \nabla f_E(p)) \in \mathbb{R}^{n+1}$ is independent of spatial location, this will prove the claim. To see that the distribution of the pair $(f(p), \nabla f_E(p))$ is independent of spatial location, note that

$$f(p) = F(y_1(p), \dots, y_k(p))$$

and

$$\nabla f_{E,i}(p) = E_i f(p) = \sum_{j=1}^k E_i y_j(p) \cdot \frac{\partial F}{\partial y_j}\Big|_{(y_1(p), \dots, y_k(p))}.$$

Since $E$ is an orthonormal set of frames, the $E_i y_j(p)$ are i.i.d. standard normal random variables, as are the $(y_1(p), \dots, y_k(p))$ and the two sets of Gaussian random variables are independent of each other.

We now turn to (2.6). Conditional on $y$ and $y_*$, $\nabla^2 f$ is a Gaussian double form. Further (cf. Section 2.4 of [18]),

$$\mathbb{E}[\nabla^2 f | y, y^*] = y^* \nabla^2 F - I\left(\sum_{i=1}^k y_i \frac{\partial F}{\partial y_i}\right)$$
$$= y^* \nabla^2 F + I \langle \nabla F, -\nabla \|x\|^2/2 \rangle |_{x=y}$$

$$= y^* \nabla^2 F + I \langle \nabla F(y), y \rangle,$$

$$\mathbb{E}[(\nabla^2 f - \mathbb{E}[\nabla^2 f | y, y_*])^2 | y, y_*] = -(I^2 + 2R) \sum_{i=1}^{k} \left( \frac{\partial F}{\partial y_i} \right)^2$$

$$= -(I^2 + 2R) \|\nabla F\|^2.$$

By Lemma 2.1, it follows that

$$\frac{1}{n!} \mathbb{E}[(-\nabla^2 f)^n | y, y^*]$$

$$= \sum_{k=0}^{\lfloor n/2 \rfloor} \frac{(-1)^{n+k}}{(n-2k)! k! 2^k} (y^* \nabla^2 F + \langle \nabla F(y), y \rangle I)^{n-2k} \|\nabla F\|^{2k} (I^2 + 2R)^k$$

$$(2.7) \quad = \sum_{j=0}^{\lfloor n/2 \rfloor} \frac{1}{j!} (-R)^j \|\nabla F\|^n$$

$$\times \sum_{l=0}^{n-2j} \frac{(y^*(-\nabla^2 F/\|\nabla F\|))^l}{l!}$$

$$\times \frac{(-1)^{n-2j-l}}{(n-2j-l)!} H_{n-2j-l}(\langle \nabla F(y)/\|\nabla F(y)\|, y \rangle) I^{n-2j-l}.$$

By Corollary 2.4,

$$\mathbb{E}[(y^* \nabla^2 F)^l | y, \nabla f_E] = \mathbb{E}[(y^* \nabla^2 F)^l | y, \langle y_* E_1, \nabla F \rangle, \ldots, \langle y_* E_n, \nabla F \rangle]$$

$$= \mathrm{Tr}^{\nabla F^\perp}(\nabla^2 F_{|\nabla F^\perp})^l I^l + \mathrm{Err}_l(\nabla f_E, y),$$

where

$$\Lambda^{l,l}(M) \ni \mathrm{Err}_l(\nabla f_E, y) = O(\|\nabla f_E\|_{\mathbb{R}^n}^2 \|(\nabla^2 F)^l\|_{\otimes^{2l} \mathbb{R}^k})$$

in the sense described in Corollary 2.4. Combining this with (2.7)

$$\frac{1}{n!} \mathbb{E}[(-\nabla^2 f)^n | y, \nabla f_E]$$

$$= \left( \sum_{j=0}^{\lfloor n/2 \rfloor} \frac{1}{j!} (-R)^j I^{n-2j} \|\nabla F\|^n \right.$$

$$\times \sum_{l=0}^{n-2j} \frac{(-1)^{n-2j-l}}{(n-2j-l)!} H_{n-2j-l}(\langle \nabla F(y), y \rangle)$$

$$\left. \times \frac{1}{l!} \mathrm{Tr}^{\nabla F^\perp}(-\nabla^2 F_{|\nabla F^\perp}/\|\nabla F\|)^l \right) + \mathrm{Err}^n(\nabla f_E, y)$$

where, up to constant factors

$$\mathrm{Err}^n(\nabla f_E, y) = \sum_{j=0}^{\lfloor n/2 \rfloor} \sum_{l=0}^{n-2j} R^j I^{n-2j-l} H_{n-2j-l}(\langle \nabla F(y), y \rangle) \mathrm{Err}_l(\nabla f_E, y),$$

from which the claim readily follows. $\square$

We have now reduced the problem of calculating expected Euler characteristics for Gaussian-related fields to the isotropic case. The following lemma gives an explicit integral representation for these EC densities.

LEMMA 2.6. *Suppose $f = F \circ z$ satisfies the conditions of Lemma 2.5 where $z = (z_1, \ldots, z_k)$ are i.i.d. isotropic, zero-mean unit-variance suitably regular Gaussian fields on $\mathbb{R}^n$ that induce the standard Riemannian structure on $\mathbb{R}^n$. Further, suppose that, for some $\delta > 0$:*

 (i) *there exist $C_1, C_2 > 0$ such that on $F^{-1}(u - \delta, u + \delta)$*

$$C_1 < \|\nabla F(x)\| < C_2;$$

 (ii) *$\nabla F$ is Lipschitz on $F^{-1}(u - \delta, u + \delta)$;*

 (iii) *the functions*

$G_{n,l,F}(z)$

$\quad = \displaystyle\int_{F^{-1}\{z\}} \|\nabla F(x)\| \widetilde{G}_{n,l,F}(x)(2\pi)^{-k/2} e^{-\|x\|^2/2} \, d\mathcal{H}_{k-1}(x)$

$\quad \triangleq \displaystyle\int_{F^{-1}\{z\}} \|\nabla F(x)\| H_{n-1-l}(\langle \nabla F(y), y \rangle(x))$

$\qquad\qquad \times \mathrm{Tr}^{\nabla F^{\perp}(x)}(-\nabla^2 F(x)/\|\nabla F(x)\|)^l (2\pi)^{-k/2} e^{-\|x\|^2/2} \, d\mathcal{H}_{k-1}(x)$

*are continuous on $(u - \delta, u + \delta)$;*

 (iv) *for all $n, l$*

$$\lim_{\varepsilon \to 0} \frac{1}{2\varepsilon} \mathbb{E}[\mathbb{1}_{\{|F(y)-u|<\varepsilon\}} |H_{n-1-l}(-\langle \nabla F(y), y \rangle)| \|\nabla^2 F(x)\|_{\otimes^{2l}\mathbb{R}^k}|] < \infty.$$

*Then,*

$$\widetilde{\rho}_n(F, u) = \frac{1}{(2\pi)^{n/2}} \sum_{l=0}^{n-1} (-1)^{n-1-l} \binom{n-1}{l} G_{n,l,F}(u).$$

PROOF. Using the original representation of the EC densities in [1] for random fields on $\mathbb{R}^n$

$$\widetilde{\rho}_n(F, u) = \rho_{f,n}(u)$$

$$= \mathbb{E}\left[\left(\frac{\partial f}{\partial t_n}\right)^+ \det(-\nabla^2 f_{|n-1})\Big| f = u, \nabla f_{|n-1} = 0\right]\varphi_{f,\nabla f_{|n-1}}(u,0)$$

$$= \lim_{\varepsilon \to 0}\frac{1}{2\varepsilon^n \omega_{n-1}}\mathbb{E}\left[\mathbb{1}_{\{|f-u|<\varepsilon\}}\mathbb{1}_{\{\|\nabla f_{|n-1}\|<\varepsilon\}}\left(\frac{\partial f}{\partial t_n}\right)^+\det(\nabla^2 f_{|n-1})\right],$$

where

$$\nabla f_{|n-1} = \left(\frac{\partial f}{\partial t_1}, \ldots, \frac{\partial f}{\partial t_{n-1}}\right)$$

is the vector made up of the first $n-1$ coordinates of $\nabla f$, $\varphi_{f,\nabla f_{|n-1}}$ is the joint density of $(f, \nabla f_{|n-1})$ and $\nabla^2 f_{|n-1}$ is the matrix of the Hessian of $f$ in the first $n-1$ coordinates. We can rewrite the expression inside the expectation above as

$$\left(\frac{\partial f}{\partial t_n}\right)^+\det(-\nabla^2 f_{|n-1})$$

$$= \frac{1}{(n-1)!}(\langle \nabla F, y_*(\partial/\partial t_n)\rangle)^+$$

$$\times (-\nabla^2 f)^{n-1}((\partial/\partial t_1, \ldots, \partial/\partial t_{n-1}), (\partial/\partial t_1, \ldots, \partial/\partial t_{n-1})).$$

From calculations similar to those in Lemma 2.5 (just restricting $\nabla^2 f$ to the first $n-1$ coordinates and reapplying Corollary 2.4), we see that

$$\mathbb{E}\left[\left(\frac{\partial f}{\partial t_n}\right)^+\det(-\nabla^2 f_{|n-1})\Big| y, \frac{\partial f}{\partial t_i}, 1 \le i \le n-1\right]$$

$$= \frac{\|\nabla F\|^n}{(2\pi)^{1/2}}\sum_{l=0}^{n-1}(-1)^{n-1-l}\binom{n-1}{l}H_{n-1-l}(\langle\nabla F(y), y\rangle)$$

$$\times \mathrm{Tr}^{\nabla F^\perp}(-\nabla^2 F/\|\nabla F\|)^l + \mathrm{Err}^{n-1}(\nabla f_{|n-1}, y)$$

$$= \frac{\|\nabla F(y)\|^n}{(2\pi)^{1/2}}\sum_{l=0}^{n-1}(-1)^{n-1-l}\binom{n-1}{l}\widetilde{G}_{n,l,F}(y) + \mathrm{Tr}(\mathrm{Err}^{n-1}(\nabla f_{|n-1}, y)),$$

where

$$\mathrm{Err}^{n-1}(\nabla f_{|n-1}, y)$$

$$= \sum_{l=0}^{n-1}I^{n-1-l}H_{n-1-l}(\langle\nabla F(y), y\rangle)O(\|\nabla f_{|n-1}\|_{\mathbb{R}^{n-1}}^2\|(\nabla^2 F)^l\|_{\otimes^{2l}\mathbb{R}^k}).$$

Assumptions (i) and (iv) above imply that

$$\lim_{\varepsilon\to 0}\frac{1}{2\varepsilon^n\omega_{n-1}}\mathbb{E}[\mathbb{1}_{\{|F(y)-u|<\varepsilon\}}\mathbb{1}_{\{\|f_{|n-1}\|<\varepsilon\}}|\mathrm{Tr}(\mathrm{Err}^{n-1}(\nabla f_{|n-1}, y))|] = 0.$$

Passing to the limit for the remaining terms

$$\widetilde{\rho}_{F,n}(u) = \lim_{\varepsilon \to 0} \frac{1}{2\varepsilon^n \omega_{n-1}}$$

$$\times \sum_{l=0}^{n-1} (-1)^{n-1-l} \binom{n-1}{l}$$

$$\times \mathbb{E}\left[\mathbb{E}\left[\mathbb{1}_{\{\|\nabla f_{|n-1}\| < \varepsilon\}} \mathbb{1}_{\{|f-u| < \varepsilon\}} \frac{\|\nabla F(y)\|^n}{(2\pi)^{1/2}} \widetilde{G}_{n,l,F}(y) \Big| y\right]\right]$$

$$= \lim_{\varepsilon \to 0} \frac{1}{2\varepsilon} \sum_{l=0}^{n-1} (-1)^{n-1-l} \binom{n-1}{l}$$

$$\times \mathbb{E}\left[\widetilde{D}(y,\varepsilon) \mathbb{1}_{\{|F(y)-u| < \varepsilon\}} \frac{\|\nabla F(y)\|}{(2\pi)^{1/2}} \widetilde{G}_{n,l,F}(y)\right],$$

where

$$\widetilde{D}(y,\varepsilon) = \frac{\gamma_{\mathbb{R}^{n-1}}\big(B_{\mathbb{R}^{n-1}}(0,\varepsilon/\|\nabla F(y)\|)\big)}{\omega_{n-1}(\varepsilon/\|\nabla F(y)\|)^{n-1}}.$$

By assumption (i),

$$\lim_{\varepsilon \to 0} \widetilde{D}(y,\varepsilon) = (2\pi)^{-(n-1)/2}$$

for all $y \in F^{-1}(u-\delta, u+\delta)$. Further for all such $y$,

$$|\widetilde{D}(y,\varepsilon) - (2\pi)^{-(n-1)/2}| < \sup_{0 < r \le \varepsilon/C_1} \left|\frac{\gamma_{\mathbb{R}^{n-1}}\big(B_{\mathbb{R}^{n-1}}(0,r)\big)}{\omega_{n-1}r^{n-1}} - (2\pi)^{-(n-1)/2}\right| \triangleq L(\varepsilon)$$

where the function $L(\varepsilon)$ is increasing and continuous on $[0, +\infty)$ and

$$\lim_{\varepsilon \to 0} L(\varepsilon) = 0.$$

Therefore, for $\varepsilon < \delta$,

$$\left| \frac{1}{2\varepsilon} \mathbb{E}[\mathbb{1}_{\{|F(y)-u| < \varepsilon\}} \widetilde{D}(y,\varepsilon) \|\nabla F(y)\| \widetilde{G}_{n,l,F}(y)] \right.$$

$$- (2\pi)^{-(n-1)/2} \frac{1}{2\varepsilon} \mathbb{E}[\mathbb{1}_{\{|F(y)-u| < \varepsilon\}} \|\nabla F(y)\| \widetilde{G}_{n,l,F}(y)] \Big|$$

$$< \frac{L(\varepsilon)}{2\varepsilon} \mathbb{E}[\mathbb{1}_{\{|F(y)-u| < \varepsilon\}} \|\nabla F(y)\| |\widetilde{G}_{n,l,F}(y)|].$$

The above implies that if, for some $\widetilde{\delta} > 0$,

$$(2.8) \qquad \sup_{0 < \varepsilon < \widetilde{\delta}} \frac{1}{2\varepsilon} \mathbb{E}[\mathbb{1}_{\{|F(y)-u| < \varepsilon\}} \|\nabla F(y)\| \widetilde{G}_{n,l,F}(y)] < +\infty,$$

we need only evaluate

$$\lim_{\varepsilon \to 0} \frac{1}{2\varepsilon} \mathbb{E}[\mathbb{1}_{\{|F(y)-u|<\varepsilon\}} \|\nabla F(y)\| \widetilde{G}_{n,l,F}(y)].$$

The Lipschitz assumption on $\nabla F$ allows us to apply Federer's coarea formula [5, 11] so that

$$
\begin{aligned}
&\lim_{\varepsilon \to 0} \frac{1}{2\varepsilon} \mathbb{E}[\mathbb{1}_{\{|F(y)-u|<\varepsilon\}} \|\nabla F(y)\| \widetilde{G}_{n,l,F}(y)] \\
&= \lim_{\varepsilon \to 0} \frac{1}{2\varepsilon} \int_{(u-\varepsilon,u+\varepsilon)} \int_{F^{-1}\{z\}} \widetilde{G}_{n,l,F}(x) e^{-\|x\|^2/2} \, d\mathcal{H}_{k-1}(x) \, dz \\
&= \lim_{\varepsilon \to 0} \frac{1}{2\varepsilon} \int_{(u-\varepsilon,u+\varepsilon)} G_{n,l,F}(z) \, dz \\
&= G_{n,l,F}(u).
\end{aligned}
$$

The Lipschitz assumption on $\nabla F$ and the coarea formula also imply (2.5). $\square$

**3. Steiner formulae for Lebesgue and other measures.** In the previous section we derived an integral representation for $\widetilde{\rho}_j(F, u)$. Ultimately, we will connect this integral representation with a certain volume of tubes expansion, to which this section is devoted. In this section we generalize the Weyl–Steiner tube formulae in order to compute non-Lebesgue volume of tubes, in particular Gaussian volumes.

We will limit ourselves to tubes around $C^2$ domains $D$ in $\mathbb{R}^k$ with outward pointing unit normal vector field $\nu$. That is, $D$ is a closed set with nonempty interior such that $\partial D$ is an embedded $C^2$ hypersurface. We write

$$T(D,r) = \{y \in \mathbb{R}^k : d(y,D) \le r\}$$

and

$$P_D(y) = \operatorname*{arg\,min}_{x \in D} d(x,y),$$

the metric projection onto $D$ when the minimum is unique. We recall the definition of the critical radius of $D$,

$$r_c(D) = \sup\{r \ge 0 : d(y,D) \le r \Rightarrow P_D(y) \text{ is unique}\},$$

and we will assume throughout that $r_c(D) > 0$.

The following lemma, whose proof we omit, summarizes some facts about the distance function of $D$, particularly the growth of the eigenvalues of the shape operator of the hypersurfaces at distance $d$ from $D$. For a hypersurface $M \subset \mathbb{R}^k$ with outward pointing normal $\nu$ and $y \in M$ the shape operator $S_y$ at $y$ is defined by

$$S_y(X_y, Z_y) = -\langle X_y, \nabla_{Z_y} \nu \rangle,$$

where $\nabla$ represents the usual covariant derivative on $\mathbb{R}^k$.

LEMMA 3.1.   *Suppose $F_D(y) = d(y, \partial D)$ for some $C^2$ domain $D$. Then,*

$$\nabla F_D(y) = \mathrm{T}_{d(y,D)\nu_{P_D(y)}}(\nu_{P_D(y)})$$

*where $\mathrm{T}_x : \mathcal{T}^*(T_y \mathbb{R}^k) \to \mathcal{T}^*(T_{x+y}\mathbb{R}^k)$ represents translation from one tangent space to another. Further, let $\lambda_1(y), \ldots, \lambda_{n-1}(y)$ be the eigenvalues of the shape operator $S_y$ acting on $T_y F^{-1}\{d(y, D)\}$; then*

$$\lambda_i(y) = \frac{\lambda_i(P_D(y))}{1 - d(y, D)\lambda_i(P_D(y))}.$$

Lemma 3.1 describes the geometry of the hypersurfaces at distance $r$ from $D$ and can be used to derive Steiner's formula, Theorem 3.2, which gives an expression for the Lebesgue volume of the tube $T(D, r)$ around $D$, in terms of the Lipschitz–Killing curvature measures $(\mathcal{L}_j(D, \cdot))_{0 \le j \le k}$ of $D$. The interested reader is referred to [11, 12, 13, 16, 19] for additional details.

THEOREM 3.2 (Steiner's formula).   *If we set*

$$A_r = \{(p, X_p) \in N(\partial D) : p \in A \cap \partial D, X_p = r\nu_p\},$$

*where $N(\partial D)$ is the normal bundle of $\partial D$ in $\mathbb{R}^k$, then Steiner's formula states that for $r < r_c(D)$*

$$(3.1) \qquad \mathcal{H}_{k-1}(\exp_{\partial D}(A_r)) = \sum_{j=1}^{n} \frac{r^{j-1}}{(j-1)!} \mathcal{M}_j^{\lambda_{\mathbb{R}^k}}(D, A).$$

In the next corollary, which is just an application of Steiner's formula, we derive a formal Taylor series expansion of the integral of smooth functions over $T(D, r)$. The functions are assumed to be in $\mathcal{S}(\mathbb{R}^k)$, the Schwartz space of functions on $\mathbb{R}^k$.

COROLLARY 3.3.   *Suppose $D$ is a $C^2$ domain in $\mathbb{R}^k$ with critical radius $r_c(D) > 0$. Further, suppose $D$ satisfies, for any $0 \le j \le k$,*

$$\int_{\partial D} \frac{1}{1 + \|x\|^\beta} \mathcal{M}_j^{\lambda_{\mathbb{R}^k}}(D, dx) < \infty,$$

*for some $\beta > 0$. Then, for any $n$, $r < \min(r_c(D), 1)$, $\varphi \in \mathcal{S}(\mathbb{R}^k)$ there exist $(\mathcal{C}_l(D; \varphi))_{l \ge 0}$ and a constant $K_n(D, \varphi)$ such that*

$$(3.2) \quad \left| \int_{T(D,r)} \varphi \, d\lambda_{\mathbb{R}^k} - \int_D \varphi \, d\lambda_{\mathbb{R}^k} - \sum_{l=1}^{n+k} \frac{r^l}{l!} \mathcal{C}_l(D; \varphi) \right| < \frac{r^{n+2}}{(n+2)!} K_n(D, \varphi).$$

*Specifically,*

$$\mathcal{C}_l(D;\varphi) = \sum_{m=0}^{l-1} \int_{\partial D} \binom{l-1}{m} \frac{d^{l-1-m}\varphi}{d\nu^{l-1-m}}\Big|_x \mathcal{M}_{m+1}^{\lambda_{\mathbb{R}^k}}(D, dx),$$

*where $d/d\nu$ represents differentiation in the direction of the outward unit normal $\nu$ and $\mathcal{M}_j^{\lambda_{\mathbb{R}^k}}(D, \cdot)$ is defined to be zero for $j > k$.*

PROOF. As $\varphi \in \mathcal{S}(\mathbb{R}^k)$ for $n > 0$, at each $p \in \partial D$, $\varphi(\exp_{\partial D}(p, r\nu_p))$ can be expressed as

$$(3.3) \qquad \varphi(\exp_{\partial D}(p, r\nu_p)) = \sum_{j=0}^{n} \frac{r^j}{j!} \frac{d^j\varphi}{d\nu^j}\Big|_p + \frac{r^{n+1}}{(n+1)!} \frac{d^{n+1}\varphi}{d\nu^{n+1}}\Big|_{\alpha(r,p)},$$

where $\alpha(r,p) = p + \theta(r,p)r\nu_p$ for some $\theta : [0, r_c(D)] \times \partial D \to (0, 1)$.

Combining (3.1) and (3.3) we see

$$\int_{\exp_{\partial D}(A_r)} \varphi(y) \, d\mathcal{H}_{k-1}(y)$$

$$= \sum_{l=0}^{n+k-1} \frac{r^l}{l!} \sum_{m=0}^{l} \binom{m}{l} \int_A \frac{d^m\varphi}{d\nu^m}\Big|_p \mathcal{M}_{l-m+1}^{\lambda_{\mathbb{R}^k}}(D, dp)$$

$$+ \sum_{l=0}^{k-1} r^l \int_A \frac{r^{n+1}}{(n+1)!} \frac{d^{n+1}\varphi}{d\nu^{n+1}}\Big|_{\alpha(r,p)} \mathcal{M}_{l+1}^{\lambda_{\mathbb{R}^k}}(D, dp)$$

and

$$\sum_{l=0}^{k-1} r^l \int_A \frac{r^{n+1}}{(n+1)!} \frac{d^{n+1}\varphi}{d\nu^{n+1}}\Big|_{\alpha(r,p)} \mathcal{M}_{l+1}^{\lambda_{\mathbb{R}^k}}(D, dp)$$

$$\leq \frac{r^{n+1}}{(n+1)!} C_k \max_{1 \leq i \leq k} \sup_x \left( (1 + \|x\|^\beta) \Big| \frac{d^{n+1}\varphi}{dx_i^{n+1}}\Big|_x \Big| \right)$$

$$\times \sum_{l=0}^{k-1} r^l \int_A \frac{1}{1 + |\alpha(r,p)|^\beta} \mathcal{M}_{l+1}^{\lambda_{\mathbb{R}^k}}(D, dp)$$

$$\leq \frac{r^{n+1}}{(n+1)!} K_n(D, \varphi).$$

Integrating over $\partial D$ and $[0, r]$ we are left with (3.2), which concludes the proof. $\square$

The following corollary gives the expression for $\mathcal{C}_l(D; \varphi)$ when $\varphi$ is the density of $\gamma_{\mathbb{R}^k}$ with respect to $\lambda_{\mathbb{R}^k}$, that is,

$$\varphi(x) = \frac{1}{(2\pi)^{k/2}} e^{-\|x\|^2/2}.$$

The $\mathcal{C}_l(D; \varphi)$ are coefficients in the promised power series expansion of the Gaussian volume of the tube $T(D, r)$.

COROLLARY 3.4. *If $\gamma_{\mathbb{R}^k}$ is the canonical Gaussian measure on $\mathbb{R}^k$, we define, for $l \geq 1$,*

$$\mathcal{M}_l^{\gamma_{\mathbb{R}^k}}(D) \triangleq \mathcal{C}_l(D; (2\pi)^{-k/2} e^{-\|x\|^2/2})$$

$$(3.4) \qquad = \frac{1}{(2\pi)^{k/2}} \sum_{m=0}^{l-1} \int_{\partial D} (-1)^m \frac{(l-1)!}{m!}$$

$$\times H_m(\langle p, \nu_p \rangle) e^{-\|p\|^2/2} \mathcal{M}_{l-m}^{\lambda_{\mathbb{R}^k}}(D, dp)$$

*and*

$$\mathcal{M}_0^{\gamma_{\mathbb{R}^k}}(D) \triangleq \gamma_{\mathbb{R}^k}(D).$$

*If, further,*

$$\int_{\partial D} \frac{1}{1 + \|p\|^\beta} \mathcal{M}_j^{\lambda_{\mathbb{R}^k}}(D, dp) < \infty,$$

*for some $\beta > 0$, then the following formal expansion holds:*

$$\gamma_{\mathbb{R}^k}(T(D, r)) = \gamma_{\mathbb{R}^k}(D) + \sum_{j=1}^{\infty} \frac{r^j}{j!} \mathcal{M}_j^{\gamma_{\mathbb{R}^k}}(D).$$

### 3.1. *Normalization of $\mathcal{M}_j^{\gamma_{\mathbb{R}^k}}(\cdot)$.*

We conclude this section by showing that the functionals $\mathcal{M}_j^{\gamma_{\mathbb{R}^k}}(\cdot)$ can be seen to be normalized independent of the dimension $k$, and can thus be extended, at least formally, to sets in $\mathbb{R}^{\mathbb{N}}$. Because of this normalization, we would be justified in writing $\mathcal{M}^\gamma$ instead of $\mathcal{M}^{\gamma_{\mathbb{R}^k}}$ in what follows, though we keep the latter notation in order to minimize confusion.

For a finite subset $\mathcal{I}$ of $\mathbb{N}$, we denote the projection from $\mathbb{R}^{\mathbb{N}}$ onto $\mathbb{R}^{\mathcal{I}}$ by $\pi_{\mathcal{I}}$. For $D \subset \mathbb{R}^k$, with $k = \#\mathcal{I}$, we can define the functionals $\mathcal{M}_j^\gamma(\cdot)$ on cylindrical sets of the form $\pi_{\mathcal{I}}^{-1}(D)$ as

$$\mathcal{M}_j^\gamma(\pi_{\mathcal{I}}^{-1}(D)) \triangleq \mathcal{M}_j^{\gamma_{\mathbb{R}^k}}(D).$$

The functionals $\mathcal{M}_j^\gamma(\cdot)$ can then hypothetically be extended by taking limits of cylindrical sets, though we will not pursue this matter here, except to show that these limits exist for at least $j = 0, 1, 2$. Of course, the functional $\mathcal{M}_0^\gamma$ is well defined, since it is just the infinite-dimensional i.i.d. Gaussian measure on $\mathbb{R}^{\mathbb{N}}$.

For $j > 0$, however, it is not immediately clear that the extension is possible. The following corollary, which is really just an application of the coarea formula in [5] and whose proof we omit, shows why such extensions are possible, at least in the cases $j = 1, 2$. By the coarea formula the expression for $\mathcal{M}_2^{\gamma_{\mathbb{R}^k}}(D)$ is the integral of what is referred to as the mean Gaussian curvature of $\partial D$ in [4], and can be extended to codimension-1 submanifolds of Wiener space and other infinite-dimensional Gaussian measure spaces.

COROLLARY 3.5. *Suppose the $C^2$ domain $D$ is given by $F^{-1}[u, +\infty)$. If both $\mathcal{M}_1^{\gamma_{\mathbb{R}^k}}(F^{-1}[u, +\infty))$ and $\mathcal{M}_2^{\gamma_{\mathbb{R}^k}}(F^{-1}[u, +\infty)) < \infty$, then*

$$\mathcal{M}_1^{\gamma_{\mathbb{R}^k}}(F^{-1}[u, +\infty)) = \mathbb{E}[\|\nabla F\| \mid F = u]\varphi_F(u),$$

$$\mathcal{M}_2^{\gamma_{\mathbb{R}^k}}(F^{-1}[u, +\infty)) = \mathbb{E}\left[-LF + \frac{\nabla^2 F(\nabla F, \nabla F)}{\|\nabla F\|^2} \Big| F = u\right]\varphi_F(u),$$

*where $\varphi_F(\cdot)$ is the density of the random variable $F(Z_1, \ldots, Z_k)$ for $Z_i$, i.i.d. N(0,1) and*

$$LF(x) = \sum_{i=1}^k \frac{\partial^2 F(x)}{\partial x_i^2} - \sum_{i=1}^k \frac{\partial F(x)}{\partial x_k}x_k$$

*is the Ornstein–Uhlenbeck operator on $C^2(\mathbb{R}^k)$.*

## 4. EC densities and Weyl–Steiner formulae: the Gaussian KFF.

In Section 2 we derived an expression for the EC densities of Gaussian-related fields involving the integral of certain functions $\widetilde{G}_{n,l,F}$ over the surface $F^{-1}\{u\}$. In Section 3 we derived expressions for coefficients in a formal series expansion of the canonical Gaussian volume of a tube around a $C^2$ domain $D$. In this section we will show that these expressions agree, up to a constant. Specifically, we show that for suitable $F \in C^2(\mathbb{R}^k)$

$$\widetilde{\rho}_j(F, u) = (2\pi)^{-j/2}\mathcal{M}_j^{\gamma_{\mathbb{R}^k}}(F^{-1}[u, +\infty))$$

which by Lemma 2.5 translates into the following Gaussian KFF:

$$(4.1) \qquad \int_\Omega \chi(M \cap y(\omega)^{-1}D)\, d\mathbb{P}(\omega) = \sum_{j=0}^n (2\pi)^{-j/2}\mathcal{L}_j(M)\mathcal{M}_j^{\gamma_{\mathbb{R}^k}}(D).$$

THEOREM 4.1 (Gaussian KFF). *Suppose $y = (y_1, \ldots, y_k)$ are i.i.d. real-valued zero-mean unit-variance suitably regular (cf. [18]) Gaussian random fields on a $C^3$ $n$-dimensional manifold $M$ that satisfy*

$$\mathbb{P}\left[\sup_{q \in B_\tau(p,h)} \|\widetilde{y}_i(p) - \widetilde{y}_i(q)\|_2 > \varepsilon\right] = o(h^n)$$

*for any $h < h_0$ and any $\varepsilon > 0$ where*

$$\widetilde{y}_i(p) = (y_i(p), \nabla y_{i,E}(p), \nabla^2 y_{i,E}(p)) \in \mathbb{R} \times \mathbb{R}^n \times^{\otimes^2} \mathbb{R}^n,$$

$$\|(x, V, H)\|_2 = |x| + \|V\|_{\mathbb{R}^n} + \|H\|_{\otimes^2 \mathbb{R}^n},$$

*and $\nabla f_E$ (resp. $\nabla^2 f_E$) are the coefficients of $\nabla f$ (resp. $\nabla^2 f$) read off in a fixed orthonormal frame $(E_1, \ldots, E_n)$.*

Let $f_D = F_D \circ y$ where $F_D$ is the distance function of a $C^2$ domain in $\mathbb{R}^k$, not necessarily compact, with $r_c(D) > 0$ which satisfies:

(i) *the functions*

$$G^{\pm}(z) = \int_{F_D^{-1}\{z\}} e^{z\langle \nu, x \rangle - \|x\|^2/2} \, d\mathcal{H}_{k-1}(x)$$

*are continuous on some neighborhood of zero where $\nu$ is the outward pointing unit normal vector field;*

(ii) *the second fundamental form $S$ of $\partial D$ is bounded, that is,*

$$|S(X, Y)_x| \le K\|X_x\|\|Y_x\|$$

*for some $K > 0$ and all $x \in \partial D$.*

*Then,*

$$
\begin{aligned}
(4.2) \quad \mathbb{E}[\chi(M \cap f_D^{-1}[0, +\infty))] &= \mathbb{E}[\chi(M \cap y^{-1}D)] \\
&= \sum_{j=0}^{n} \mathcal{L}_j(M) \frac{1}{(2\pi)^{j/2}} \mathcal{M}_j^{\gamma_{\mathbb{R}^k}}(D).
\end{aligned}
$$

*In particular, the above relation holds for every compact $C^2$ domain $D$, and its complement $D^c$.*

REMARK. The conditions above are not necessary; indeed, there are cases such as the $F$ field defined in [20] (cf. Section 5.3) where the boundedness condition on the second fundamental form of the boundary of the domain is not satisfied but the EC densities exist and are given by

$$\sum_{l=0}^{n-1} (-1)^{n-1-l} \binom{n-1}{l} \lim_{\varepsilon \to 0} \frac{1}{2\varepsilon} \mathbb{E}[\mathbb{1}_{\{|F(y)-u|\}} \|\nabla F\| G_{n,l,F}(y)].$$

However, it should be noted that in this case, the above coefficients are not the true coefficients in the corresponding power series, because the domain $D_u$, defined in Section 5.3, has $r_c(D_u) = 0$.

REMARK. Using the generalized Morse theorem of [16] or stratified Morse theory [3], Theorem 4.1 can be extended to include piecewise smooth domains and/or submanifolds of $\mathbb{R}^k$ (i.e., the domains in the space where the random fields take values). We give an example of a piecewise smooth convex domain $D$ in Section 5.4, where we calculate the $\mathcal{M}_j^{\gamma_{\mathbb{R}^k}}(D)$ for $D$ a cone in $\mathbb{R}^2$ and show that the result agrees with known results for right-angled cones [24]. One approach to these generalizations, following the notions of continuity of the Lipschitz–Killing curvature measures as in [11], would be to construct a limiting argument to justify the following computation for certain smooth $D$ sets with positive critical radius in $\mathbb{R}^k$:

$$\int_\Omega \chi(M \cap y(\omega)^{-1}D)\,d\mathbb{P}(\omega)$$

$$= \lim_{r \to 0} \int_\Omega \chi(M \cap y(\omega)^{-1}T(D,r))\,d\mathbb{P}(\omega)$$

$$= \lim_{r \to 0} \sum_{i=0}^n \mathcal{L}_i(M)\frac{1}{(2\pi)^{i/2}}\mathcal{M}_i^{\gamma_{\mathbb{R}^k}}(T(D,r)) = \sum_{i=0}^n \mathcal{L}_i(M)\frac{1}{(2\pi)^{i/2}}\mathcal{M}_i^{\gamma_{\mathbb{R}^k}}(D).$$

In the interest of brevity, we will not pursue these generalizations here. The above theorem, however, does specify the functional form that all of these generalizations should have, that is, the contribution of the parameter space is in the form of the Lipschitz–Killing curvature measures, which can be defined for piecewise $C^2$ submanifolds of $C^3$ manifolds, and the contribution from the Gaussian space is in the form of the coefficients in an expansion for the volume of a Gaussian tube, which can similarly be defined for piecewise $C^2$ manifolds in $\mathbb{R}^k$. For a more geometric approach to the above problem, see [3] which shows that (4.2) can be thought of as a limit of the classical KFF.

REMARK. Note that the conditions on the $y_i$'s are not overly restrictive (cf. [18]), as $C^2$ fields whose second derivatives have a covariance function satisfying the "usual" $1/(-\log(h))^{1+\delta}$ moduli of continuity conditions are included above.

PROOF. Most of the work has been done in the preliminary lemmas. The conclusion

$$\widetilde{\rho}_j(F_D, u) = (2\pi)^{-j/2}\mathcal{M}_j^{\gamma_{\mathbb{R}^k}}(D)$$

follows from the fact that $f_D$ satisfies the conditions of Lemmas 2.5 and 2.6. Lemma 3.1 ensures that the second fundamental form, or shape operator $S$ of $\partial D$ is bounded. Verifying the equality above then follows from comparing the definitions in Corollary 3.4 of $\mathcal{M}_j^{\gamma_{\mathbb{R}^k}}(D)$ and those of $G_{n,l,F_D}(0)$ in

Lemma 2.6, noting that for $x \in \partial D$ and all $X_x, Y_x \in T_x \partial D$,

$$\nabla F_D(x) = \nu_x,$$

$$\nabla^2 F_D(X_x, Y_x) = -S(X_x, Y_x).$$

The fact that the conditions of Lemma 2.5 are satisfied follows from Lemmas 2.4 and 2.5 of [18] combined with the growth conditions on $y_i$ and its derivatives as well as the boundedness of $S$. As for the conditions of Lemma 2.6, since $F_D$ is a distance function of a set with positive critical radius, conditions (ii) and (iv) (cf. [11]) of Lemma 2.6 are automatically satisfied. A straightforward calculation shows that assumptions (i) and (ii) above imply conditions (ii) and (iii) of Lemma 2.6.  □

**5. Examples.** Throughout this section, $y = (y_1, \ldots, y_k) : M \to \mathbb{R}^k$ will denote a generic sequence of i.i.d. zero-mean unit-variance Gaussian fields on $M$ satisfying the conditions of Theorem 4.1, whose dimension $k$ will vary as needed in each example.

5.1. *Real-valued Gaussian processes.* Given a unit vector $z \in \mathbb{R}^k$, we define the function $\psi_z(x) = \langle x, z \rangle$, so that $\psi_z \circ y$ is a real-valued, centered unit-variance Gaussian random field. As mentioned in Section 2, the EC densities of $\psi_z \circ y$ are given, for $j \geq 1$, by

$$(5.1) \qquad \widetilde{\rho}_j(\psi_z, u) = \frac{1}{(2\pi)^{(j+1)/2}} H_{j-1}(u) e^{-u^2/2}.$$

This result can easily be rederived in light of Theorem 4.1 as follows. Let $D_u = \psi_z^{-1}([u, +\infty))$ be a half-space in $\mathbb{R}^k$. Clearly

$$T(D_u, r) = \psi_z^{-1}[u - r, +\infty) = D_{u-r},$$

that is, a tube around a half-space in $\mathbb{R}^k$ is another half-space. As

$$\frac{d^j}{dx^j} e^{-x^2/2} = (-1)^j H_j(x) e^{-x^2/2},$$

it follows that

$$\gamma_{\mathbb{R}^k}(T(D_u, r)) = 1 - \Phi(u - r)$$

$$= 1 - \left( \Phi(u) + \sum_{j=1}^{\infty} \frac{(-r)^j}{j!} \frac{(-1)^{j-1}}{\sqrt{2\pi}} H_{j-1}(u) e^{-u^2/2} \right)$$

$$= 1 - \Phi(u) + \sum_{j=1}^{\infty} \frac{r^j}{j!} \frac{1}{\sqrt{2\pi}} H_{j-1}(u) e^{-u^2/2},$$

so that, for $j \geq 1$,

$$\mathcal{M}_j^{\gamma_{\mathbb{R}^k}}(D_u) = \frac{1}{\sqrt{2\pi}} H_{j-1}(u) e^{-u^2/2},$$

which, by Theorem 4.1, implies (5.1).

5.2. *The $\chi^2$ and noncentral $\chi^2$ case.* In this subsection we derive the EC densities of a $\chi_k^2$ random field, as defined in [1, 20]. Note that the EC densities were also derived in [20], but we rederive them here as a simple application of Theorem 4.1. We set $G(x) = \|x\|^2$; then $g = G \circ y$ is a $\chi_k^2$ random field on $M$. Note that the EC densities of $\sqrt{g}$ are related to the EC densities of $g$ by

$$\rho_{\sqrt{g},j}(\sqrt{u}) = \widetilde{\rho}_j(\sqrt{G}, \sqrt{u}) = \widetilde{\rho}_j(G, u) = \rho_{g,j}(u),$$

so that it suffices to calculate the EC densities of $\sqrt{g}$.

We set

$$D_x = \sqrt{G}^{-1}[x, +\infty) = \overline{\mathbb{R}^k \backslash B_{\mathbb{R}^k}(0, x)},$$

so that $T(D_x, r) = D_{x-r}$. Therefore,

$$\gamma_{\mathbb{R}^k}(T(D_x, r)) = \gamma_{\mathbb{R}^k}(D_{x-r}).$$

It remains to express the right-hand side above as a Taylor series in $r$.

The density $f_k = d(\sqrt{G}_*(\gamma_{\mathbb{R}^k}))/d\lambda_{\mathbb{R}}$ of the square root of a $\chi_k^2$ random variable is

$$f_k(x) = \frac{1}{\Gamma(k/2)2^{(k-2)/2}} x^{k-1} e^{-x^2/2}$$

and

$$\gamma_{\mathbb{R}^k}(D_{x-r}) = \gamma_{\mathbb{R}^k}(D_x) - \sum_{j=1}^{\infty} \frac{r^j}{j!} (-1)^j \frac{d^{j-1}f_k}{dt^{j-1}}\Big|_{t=x}.$$

Direct calculations show that

$$\frac{d^{j-1}f_k}{dt^{j-1}} = \frac{t^{k-j}e^{-t^2/2}}{\Gamma(k/2)2^{(k-2)/2}} \sum_{l=0}^{\lfloor (j-1)/2 \rfloor} \sum_{m=0}^{j-1-2l} \mathbb{1}_{\{k \geq j-m-2l\}} \binom{k-1}{j-1-m-2l}$$
$$\times \frac{(-1)^{m+l}(j-1)!}{m!l!2^l} t^{2m+2l}.$$

Combining the two equations gives

$$\gamma_{\mathbb{R}^k}(D_{x-r}) = \gamma_{\mathbb{R}^k}(D_x) + \sum_{j=1}^{\infty} \frac{r^j}{j!} \frac{x^{k-j}e^{-x^2/2}}{\Gamma(k/2)2^{(k-2)/2}}$$
$$\times \left( \sum_{l=0}^{\lfloor (j-1)/2 \rfloor} \sum_{m=0}^{j-1-2l} \mathbb{1}_{\{k \geq j-m-2l\}} \binom{k-1}{j-1-m-2l} \right.$$
$$\left. \times \frac{(-1)^{j-1+m+l}(j-1)!}{m!l!2^l} x^{2m+2l} \right).$$

We thus conclude, by Theorem 4.1, that for $j \geq 1$,

$$\widetilde{\rho}_j(\sqrt{G}, x) = \frac{x^{k-j} e^{-x^2/2}}{(2\pi)^{j/2} \Gamma(k/2) 2^{(k-2)/2}}$$

$$\times \left( \sum_{l=0}^{\lfloor (j-1)/2 \rfloor} \sum_{m=0}^{j-1-2l} \mathbb{1}_{\{k \geq j-m-2l\}} \binom{k-1}{j-1-m-2l} \right.$$

$$\left. \times \frac{(-1)^{j-1+m+l}(j-1)!}{m! l! 2^l} x^{2m+2l} \right),$$

and thus the EC densities for $j \geq 1$ of the $\chi_k^2$ random field $g$ are given by

$$\widetilde{\rho}_j(G, u) = \frac{u^{(k-j)/2} e^{-u/2}}{(2\pi)^{j/2} \Gamma(k/2) 2^{(k-2)/2}}$$

$$\times \left( \sum_{l=0}^{\lfloor (j-1)/2 \rfloor} \sum_{m=0}^{j-1-2l} \mathbb{1}_{\{k \geq j-m-2l\}} \binom{k-1}{j-1-m-2l} \right.$$

$$\left. \times \frac{(-1)^{j-1+m+l}(j-1)!}{m! l! 2^l} u^{m+l} \right),$$

which agrees with [20].

Using the above formula for the $\chi_k^2$ EC densities, we can derive the EC densities of a noncentral $\chi^2$ which we define to be a Gaussian-related field with the function $G_\mu : \mathbb{R}^k \to \mathbb{R}$

$$G_\mu(y) = \|y - \mu\|^2.$$

Recalling that the density $f_{\alpha,k}$ of the square root of a noncentral $\chi_k^2$ random variable with noncentrality parameter $\alpha$ can be expressed as

$$f_{\alpha,k}(x) = \sum_{j=0}^{\infty} e^{-\alpha/2} \frac{\alpha^j}{2^j j!} f_{k+2j}(x),$$

where $f_k(x)$ is as above, the density of a $\chi_k^2$ random variable, and noting that, in our case $\alpha = \|\mu\|^2$, we obtain the following.

LEMMA 5.1.

$$\widetilde{\rho}_j(G_\mu, u) = \sum_{i=0}^{\infty} e^{-\|\mu\|^2/2} \frac{\|\mu\|^{2i}}{2^i i!} \frac{u^{(k+2i-j)/2} e^{-u/2}}{(2\pi)^{j/2} \Gamma((k+2i)/2) 2^{(k+2i-2)/2}}$$

$$\times \left( \sum_{l=0}^{\lfloor (j-1)/2 \rfloor} \sum_{m=0}^{j-1-2l} \mathbb{1}_{\{k \geq j-m-2l-2i\}} \binom{k+2i-1}{j-1-m-2l} \right.$$

$$\left. \times \frac{(-1)^{j-1+m+l}(j-1)!}{m! l! 2^l} u^{m+l} \right).$$

5.3. *The F case.* In this subsection we derive the EC densities of an $F_{k_1,k_2}$ random field, as defined in [20]. Unlike the previous two subsections, instead of using the representation of the EC densities as coefficients in a Taylor series expansion of the volume of certain tubular neighborhoods, for the $F_{k_1,k_2}$ field we use the expression for $\mathcal{M}_j^{\gamma_{\mathbb{R}^{k_1+k_2}}}(\cdot)$ given in Corollary 3.4.

We set

$$F(y) = \frac{k_2}{k_1} \frac{\sum_{i=1}^{k_1}(y_i)^2}{\sum_{i=1}^{k_2}(y_{k_1+i})^2},$$

and define the $F_{k_1,k_2}$ random fields $f$ by $f = F \circ y$. Setting $D_u = F^{-1}[u, +\infty)$ we see

$$\partial D_u = F^{-1}\{u\} = \bigcup_{r \in \mathbb{R}^+} S_{rk_1u/k_2}(\mathbb{R}^{k_1}) \times S_r(\mathbb{R}^{k_2}).$$

The definition of the functionals $(\mathcal{M}_j^{\gamma_{\mathbb{R}^{k_1+k_2}}}(\cdot))_{j \geq 0}$, along with the fact $\langle \nabla F(y), y \rangle = 0$ for all $y \in \mathbb{R}^k$ and

$$H_n(0) = \begin{cases} 0, & n \text{ odd}, \\ (-1)^l \dfrac{(2l)!}{l!2^l}, & n = 2l, \end{cases}$$

imply that

$$
\begin{aligned}
(2\pi)^{j/2}\widetilde{\rho}_j(F,u) &= \mathcal{M}_j^{\gamma_{\mathbb{R}^{k_1+k_2}}}(D_u) \\
&= (j-1)! \sum_{l=0}^{\lfloor (j-1)/2 \rfloor} \frac{(-1)^l}{l!2^l} \frac{1}{(2\pi)^{(k_1+k_2)/2}} \\
&\qquad\qquad \times \int_{\partial D_u} e^{-\|x\|^2/2} \mathcal{M}_{j-2l}(D_u; dx) \\
&= (j-1)! \sum_{l=0}^{\lfloor (j-1)/2 \rfloor} \frac{(-1)^l}{l!2^l} \frac{1}{(2\pi)^{(k_1+k_2)/2}} \\
&\qquad\qquad \times \int_{\partial D_u} e^{-\|x\|^2/2} \frac{1}{(j-2l-1)!} \\
&\qquad\qquad\qquad \times \mathrm{Tr}^{\partial D_u}(S_{\partial D_u}^{j-2l-1}) \, d\mathcal{H}_{k_1+k_2-1}(x).
\end{aligned}
$$

We now proceed to evaluate, for $0 \leq m \leq k_1 + k_2 - 1$,

$$\mathrm{Tr}^{\partial D_F}(S_{\partial D_F}^m)(y) = \mathrm{Tr}^{\partial D_{F(y)}}((-\|\nabla F(y)\|^{-1}\nabla^2 F(y)_{|\partial D_{F(y)}})^m)$$

in terms of $U(y) = \sum_{j=1}^{k_1}(y_j)^2$ and $V(y) = \sum_{j=1}^{k_2}(y_{j+k_1})^2$ and $G(y) = U(y)/V(y)$.

LEMMA 5.2. *We have the following expression for* $\mathrm{Tr}^{\partial D_F}(S^m_{\partial D_F})$:

$$\frac{1}{m!}\mathrm{Tr}^{\partial D_F}(S^m_{\partial D_F}) = \frac{1}{m!}\mathrm{Tr}^{\partial D_F}((-\|\nabla F\|^{-1}\nabla^2 F_{|\partial D_F})^m)$$

$$= \left(\frac{1}{\sqrt{VG(1+G)}}\right)^m \sum_{i=0}^m (-1)^{m-i} G^i \binom{k_1-1}{m-i}\binom{k_2-1}{i}.$$

*Further,*

$$\frac{1}{(2\pi)^{(k_1+k_2)/2}}\int_{\partial D_u}\frac{1}{m!}\mathrm{Tr}^{\partial D_u}(S^m_{\partial D_u})e^{-\|x\|^2}\,d\mathcal{H}_{k_1+k_2-1}(x)$$

$$= \frac{\Gamma((k_1+k_2-m-1)/2)}{2^{(m-1)/2}\Gamma(k_1/2)\Gamma(k_2/2)}\left(\frac{k_1 u}{k_2}\right)^{(k_1-1-m)/2}\left(1+\frac{k_1 u}{k_2}\right)^{-(k_1+k_2-2)/2}$$

$$\times \sum_{i=0}^m (-1)^{m-i}\left(\frac{k_1 u}{k_2}\right)^i \binom{k_1-1}{m-i}\binom{k_2-1}{i}.$$

PROOF. Since $F = k_2 U/k_1 V$, a straightforward calculation shows

$$\|\nabla F\| = \frac{2k_2}{k_1}\sqrt{G(1+G)}\cdot\frac{1}{\sqrt{V}},$$

$$\|\nabla F\|^{-1}\nabla F = \sqrt{\frac{1}{1+G}}\frac{1}{2\sqrt{U}}\nabla U - \sqrt{\frac{G}{1+G}}\frac{1}{2\sqrt{V}}\nabla V,$$

$$\|\nabla F\|^{-1}\nabla^2 F = \frac{1}{2\sqrt{G(1+G)}}$$

$$\times\left(\frac{1}{\sqrt{V}}\nabla^2 U - \frac{1}{V^{3/2}}(dV\otimes dU + dU\otimes dV)\right.$$

$$\left. + \frac{2G}{V^{3/2}}\,dV\otimes dV - \frac{G}{\sqrt{V}}\nabla^2 V\right).$$

We now evaluate the matrix of $\|\nabla F\|^{-1}\nabla^2 F$ in a specific set of orthonormal frames of $\nabla F^\perp$. Considering $\mathbb{R}^{k_1+k_2}$ to be the product $\mathbb{R}^{k_1}\times\mathbb{R}^{k_2}$ with the product metric, we fix subspaces $L_1$ the kernel of $\nabla U$ in $\mathbb{R}^{k_1}$ and similarly $L_2$, the kernel of $\nabla V$ in $\mathbb{R}^{k_2}$, for which we choose orthonormal frames $B_1, B_2$. The desired set of frames for the kernel of $\nabla F$ is then $B = \{B_1, B_2, X\}$, where

$$X = \sqrt{\frac{G}{1+G}}\frac{1}{2\sqrt{U}}\nabla U - \sqrt{\frac{1}{1+G}}\frac{1}{2\sqrt{V}}\nabla V.$$

The matrix of $-\|\nabla F\|^{-1}\nabla^2 F$ in this set of frames is diagonal, with entries

$$\frac{1}{\sqrt{VG(1+G)}}(\overbrace{-1,\ldots,-1}^{k_1-1\mathrm{times}},\overbrace{G,\ldots,G}^{k_2-1\mathrm{times}},0),$$

from which the result follows by expanding the trace of the $m$th power of such a diagonal matrix.

Federer's coarea formula (cf. [5, 11]) implies

$$\frac{1}{(2\pi)^{(k_1+k_2)/2}} \int_{\partial D_u} e^{-\|x\|^2/2} \frac{1}{m!} \operatorname{Tr}^{\partial D_u}(S_{\partial D_u}^m) \, d\mathcal{H}_{k_1+k_2-1}(x)$$

$$= \lim_{\varepsilon \to 0} \frac{1}{2\varepsilon} \int_{\mathbb{R}^{k_1+k_2}} \|\nabla F\| \mathbb{1}_{\{|F-u|<\varepsilon\}} \frac{1}{m!} \operatorname{Tr}^{\partial D_F}(S_{\partial D_F}^m) \, d\gamma_{\mathbb{R}^{k_1+k_2}}(x)$$

$$= \mathbb{E}\left[\|\nabla F\| \frac{1}{m!} \operatorname{Tr}^{\partial D_F}(S_{\partial D_F}^m) \Big| F=u\right] \frac{dF_*(\gamma_{\mathbb{R}^{k_1+k_2}})}{d\lambda_{\mathbb{R}}}(u).$$

The second conclusion now follows after substituting in the density $dF_*(\gamma_{\mathbb{R}^{k_1+k_2}})/d\lambda_{\mathbb{R}}$ and noting that $V(1+G) = U+V \sim \chi^2_{k_1+k_2}$ is independent of $G$ and

$$\mathbb{E}[(U+V)^p] = 2^p \frac{\Gamma((k_1+k_2)/2 + p)}{\Gamma((k_1+k_2)/2)}. \qquad \square$$

Combining (5.2) and the previous lemma we have, in agreement with [20], the following

LEMMA 5.3.    *The EC densities for the $F_{k_1,k_2}$ random field are given by*

$$(2\pi)^{j/2} \widetilde{\rho}_j(F, u)$$

$$= \frac{\Gamma((k_1+k_2-j)/2)}{2^{(j-2)/2} \Gamma(k_1/2) \Gamma(k_2/2)} \left(\frac{k_1 u}{k_2}\right)^{(k_1-j)/2} \left(1 + \frac{k_1 u}{k_2}\right)^{-(k_1+k_2-2)/2}$$

$$\times (-1)^{j-1}(j-1)! \sum_{l=0}^{\lfloor (j-1)/2 \rfloor} \frac{\Gamma((k_1+k_2-j)/2 + l)}{\Gamma((k_1+k_2-j)/2)l!}$$

$$\times \sum_{i=0}^{j-2l-1} (-1)^{i+l} \left(\frac{k_1 u}{k_2}\right)^{i+l}$$

$$\times \binom{k_1-1}{j-1-2l-i} \binom{k_2-1}{i}.$$

5.4.  *Cones in $\mathbb{R}^2$ and correlated conjunctions.*    In this section we study conjunctions of correlated Gaussian fields, which in terms of fields can be defined in terms of the minimum of two correlated Gaussian processes. Specifically, as in the Introduction, given $y = (y_1, y_2)$ i.i.d. centered unit-variance

Gaussian processes on some manifold $M$, we form two new Gaussian processes as follows:

$$z_1 = y_1,$$
$$z_2 = \rho \cdot y_1 + \sqrt{1-\rho^2} \cdot y_2,$$

and define $\tilde{z}_1$ and $\tilde{z}_2$, our isotropic versions of these processes, as in Section 3. We define the conjunction of $z_1$ and $z_2$ at a point $p \in M$ by

$$z_1 \wedge z_2(p) \triangleq \min(z_1(p), z_2(p)).$$

It is easy to see that

$$z_1 \wedge z_2^{-1}[u, +\infty) = z_1^{-1}[u, +\infty) \cap z_2^{-1}[u, +\infty)$$
$$= y^{-1}(K(u, \rho)),$$

where $K(u, \rho)$ is a cone in $\mathbb{R}^2$ so that, to calculate the EC densities for this process, it suffices to calculate the expected Euler characteristic of $M \cap y^{-1}K$ for a general cone with arbitrary vertex in $\mathbb{R}^2$, which we proceed to do.

Define the cone

$$C(v_1, v_2, w) = \{z \in \mathbb{R}^2 : z = w + a_1 v_1 + a_2 v_2, a_1, a_2 \geq 0\}$$
$$= w + \{z \in T_w \mathbb{R}^2 : z = a_1 v_1 + a_2 v_2, a_1, a_2 \geq 0\},$$

where $T_w \mathbb{R}^2$ is the tangent space to $\mathbb{R}^2$ at $w$. A sketch of the cone appears in Figure 1. We derive an expression for the quantity

$$\mathcal{M}^\gamma(C(v_1, v_2, w))$$

which, by Theorem 4.1 and the remarks following it, allows us to compute

$$\mathbb{E}[\chi(M \cap y^{-1}C(v_1, v_2, w))].$$

Associated to $C(v_1, v_2, w)$ is its normal cone

$$C^\perp(v_1, v_2, w) = \{z \in T_w \mathbb{R}^2 : \langle z, v_1 \rangle < 0, \langle z, v_2 \rangle < 0\},$$

with link

$$L(C^\perp(v_1, v_2, w)) = \{z \in S(T_w \mathbb{R}^2) : \langle z, v_1 \rangle < 0, \langle z, v_2 \rangle < 0\},$$

where $S(T_w \mathbb{R}^2)$ is the unit sphere in $T_w \mathbb{R}^2$.

In this case, the domain $C(v_1, v_2, w)$ is not smooth; however, for any $\delta > 0$, $T(C(v_1, v_2, w), \delta)$ is at least $C^1$ and the coefficients of

$$\gamma_{\mathbb{R}^2}(T(T(C(v_1, v_2, w), \delta), r)),$$

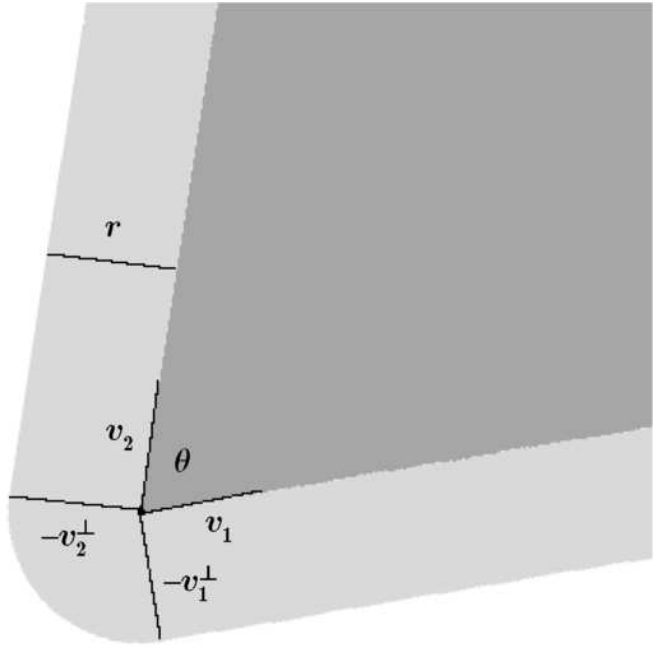

- **(0, 0)**

FIG. 1.  *Tube of radius $r$ around $C(v_1, v_2, w)$.*

as a power series in $r$ will, for small $\delta > 0$, be close to those of $\gamma_{\mathbb{R}^2}(T(C(v_1, v_2, w), r))$, that is, for small $\delta > 0$,

$$\mathbb{E}[\chi(M \cap y^{-1}C(v_1, v_2, w))] \simeq \mathbb{E}[\chi(M \cap y^{-1}(T(C(v_1, v_2, w), \delta)))]$$

$$= \frac{1}{(2\pi)^{n/2}} \sum_{j=0}^{n} \mathcal{L}_j(M) \mathcal{M}_j^{\gamma_{\mathbb{R}^2}}(T(C(v_1, v_2, w), \delta))$$

$$\simeq \frac{1}{(2\pi)^{n/2}} \sum_{j=0}^{n} \mathcal{L}_j(M) \mathcal{M}_j^{\gamma_{\mathbb{R}^2}}(C(v_1, v_2, w)),$$

where, for the cone $C(v_1, v_2, w)$, $\mathcal{M}_j^{\gamma_{\mathbb{R}^2}}(C(v_1, v_2, w))$ is defined as the coefficient of $r^j/j!$ in a Taylor series expansion of $\gamma_{\mathbb{R}^2}(T(D, r))$. The expansion is split up into integrals over three regions, depicted in Figure 2.

Without loss of generality, we can choose a basis of $T_w\mathbb{R}^2$ so that the coefficients of $v_1$ are $(1, 0)$ and those of $v_2$ are $(\cos\theta, \sin\theta)$ where $\theta = \cos^{-1}(\langle v_1, v_2 \rangle)$. We then set $v_1^\perp$ to be the unit vector orthogonal to $v_1$ such that $\langle v_1^\perp, v_2 \rangle > 0$,

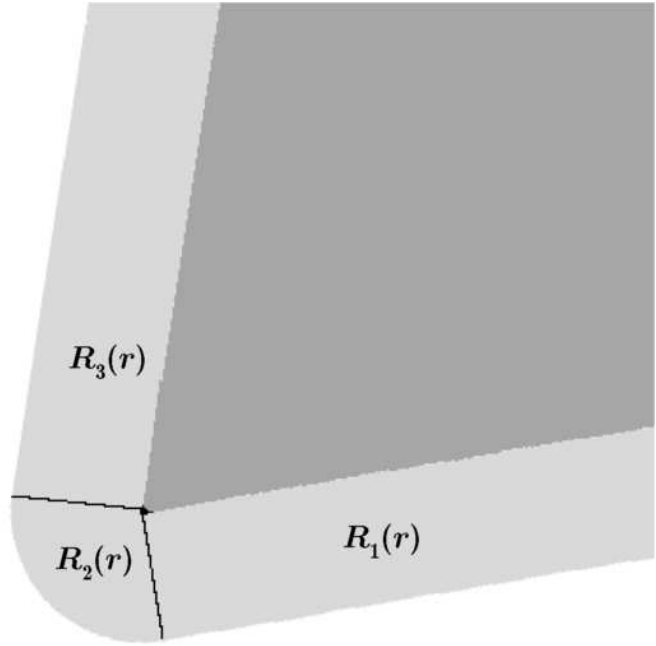

• $(0, 0)$

Fig. 2.  *Regions of integration for power series expansion of $\gamma_{\mathbb{R}^2}(T(C(v_1, v_2, w)))$.*

that is, the coefficients of $v_1^\perp$ with respect to our chosen basis are $(0, 1)$. The volume of the first region, $R_1(r)$, is thus

$$\gamma_{\mathbb{R}^2}(R_1(r)) = \gamma_{\mathbb{R}^2}(\{z \in \mathbb{R}^2 : \langle v_1^\perp, z \rangle \in [\langle v_1^\perp, w \rangle - r, \langle v_1^\perp, w \rangle], \langle v_1, z \rangle \geq \langle v_1, w \rangle\})$$

$$= \frac{(1 - \Phi(\langle v_1, w \rangle))}{\sqrt{2\pi}} \left( \sum_{j=1}^\infty \frac{r^j}{j!} H_{j-1}(\langle v_1^\perp, w \rangle) e^{-\langle v_1^\perp, w \rangle^2 / 2} \right).$$

By symmetry,

$$\gamma_{\mathbb{R}^2}(R_3(r)) = \gamma_{\mathbb{R}^2}(\{z \in \mathbb{R}^2 : \langle v_2^\perp, z \rangle \in [\langle v_2^\perp, w \rangle - r, \langle v_2^\perp, w \rangle], \langle v_2, z \rangle \geq \langle v_2, w \rangle\})$$

$$= \frac{(1 - \Phi(\langle v_2, w \rangle))}{\sqrt{2\pi}} \left( \sum_{j=1}^\infty \frac{r^j}{j!} H_{j-1}(\langle v_2^\perp, w \rangle) e^{-\langle v_2^\perp, w \rangle^2 / 2} \right).$$

Finally, using polar coordinates

$$\gamma_{\mathbb{R}^2}(R_2(r)) = \frac{1}{2\pi} \sum_{j=0}^{\infty} \int_{[0,r] \times L(C^\perp(v_1,v_2,w))} \frac{t^{j+1}}{j!} \frac{d^j e^{-\|z\|^2/2}}{d\nu^j} \Big|_{z=w} dt\, d\nu$$

$$= \frac{1}{2\pi} \sum_{j=2}^{\infty} \frac{r^j}{j!} (j-1) \int_{L(C^\perp(v_1,v_2,w))} \frac{d^{j-2} e^{-\|z\|^2/2}}{d\nu^{j-2}} \Big|_{z=w} d\nu.$$

Writing $\nu = -\sin\widetilde\theta \cdot v_1 - \cos\widetilde\theta \cdot v_1^\perp$, we have

$$\frac{d^{j-2} e^{-\|z\|^2/2}}{d\nu^{j-2}} \Big|_{z=w} = (-1)^{j-2} \left( \sin\widetilde\theta \frac{d}{dv_1} + \cos\widetilde\theta \frac{d}{dv_1^\perp} \right)^{j-2} e^{-\|z\|^2/2} \Big|_{z=w}$$

$$= \sum_{l=0}^{j-2} \binom{j-2}{l} \sin^{j-2-l}\widetilde\theta \cos^l\widetilde\theta$$

$$\times H_{j-2-l}(\langle v_1, w\rangle) H_l(\langle v_1^\perp, w\rangle) e^{-\|w\|^2/2}.$$

Noting that $\nu \in L(C^\perp(v_1,v_2,w))$ if and only if $\widetilde\theta \in (0, \pi-\theta)$ we see

$$\gamma_{\mathbb{R}^2}(R_2(r)) = \frac{1}{2\pi} \sum_{j=2}^{\infty} \frac{r^j}{j!} \sum_{l=0}^{j-2} (j-1) \binom{j-2}{l}$$

$$\times H_{j-2-l}(\langle v_1, w\rangle) H_l(\langle v_1^\perp, w\rangle) e^{-\|w\|^2/2}$$

$$\times \int_0^{\pi-\theta} \sin^{j-2-l}\widetilde\theta \cos^l\widetilde\theta\, d\widetilde\theta.$$

Straightforward calculations show

$$K_{j,l}(\theta) \triangleq (j-1) \int_0^{\pi-\theta} \sin^{j-2-l}\widetilde\theta \cos^l\widetilde\theta\, d\widetilde\theta$$

$$= \begin{cases} \dfrac{j-1}{2} IB_{(j-1-l)/2,(l+1)/2}(\sqrt{\sin\theta}), & \text{if } \theta \geq \pi/2, \\[2mm] (-1)^l \dfrac{j-1}{2}(B_{(j-1-l)/2,(l+1)/2} - IB_{(j-1-l)/2,(l+1)/2}(\sqrt{\sin\theta})) \\[1mm] \quad + \dfrac{j-1}{2} B_{(j-1-l)/2,(l+1)/2}, & \text{if } \theta < \pi/2, \end{cases}$$

where

$$IB_{\nu_1,\nu_2}(x) = \int_0^x t^{\nu_1-1}(1-t)^{\nu_2-1}\, dt$$

is the incomplete beta function and $B_{\nu_1,\nu_2} = IB_{\nu_1,\nu_2}(1)$ is the beta function.

Putting the above together, we have proved the following.

LEMMA 5.4.   *The coefficient of $r^j/j!$ in the power series expansion of*

$$\gamma_{\mathbb{R}^2}(T(C(v_1, v_2, w), r))$$

*is given by*

$$\mathcal{M}_j^{\gamma_{\mathbb{R}^2}}(C(v_1, v_2, w))$$

$$= \begin{cases}
\gamma_{\mathbb{R}^2}(C(v_1, v_2, w)), & j = 0, \\[2mm]
\dfrac{(1 - \Phi(\langle v_1, w \rangle))}{\sqrt{2\pi}} e^{-\langle v_1^\perp, w \rangle^2/2} + \dfrac{(1 - \Phi(\langle v_2, w \rangle))}{\sqrt{2\pi}} e^{-\langle v_2^\perp, w \rangle^2/2}, & \\
& j = 1, \\[2mm]
\dfrac{(1 - \Phi(\langle v_1, w \rangle))}{\sqrt{2\pi}} H_{j-1}(\langle v_1^\perp, w \rangle) e^{-\langle v_1^\perp, w \rangle^2/2} & \\
\quad + \dfrac{(1 - \Phi(\langle v_2, w \rangle))}{\sqrt{2\pi}} H_{j-1}(\langle v_2^\perp, w \rangle) e^{-\langle v_2^\perp, w \rangle^2/2} & \\
\quad + \dfrac{1}{2\pi} \displaystyle\sum_{l=0}^{j-2} \binom{j-2}{l} K_{j,l}(\theta) H_{j-2-l}(\langle v_1, w \rangle) H_l(\langle v_1^\perp, w \rangle) e^{-\|w\|^2/2}, & \\
& j \geq 2.
\end{cases}$$

It remains only to relate the above lemma to our original goal, that is, the EC densities of the field $z_1 \wedge z_2$, which amounts to determining $v_1, v_2$ and $w$ for the cone $K(u, \rho)$. We can take $v_1^\perp = (1, 0)$ and $v_2^\perp = (\rho, \sqrt{1 - \rho^2})$ so that $v_1 = (0, 1)$ and $v_2 = (\sqrt{1 - \rho^2}, -\rho)$ and $w = (u, u/\sqrt{1 + \rho})$.

**Acknowledgments.**   The author would like to thank Drs. Robert Adler and Keith Worsley for their help and encouragement as thesis advisors, which is where this work originated.

DEPARTMENT OF STATISTICS
SEQUOIA HALL
STANFORD UNIVERSITY
STANFORD, CALIFORNIA 94305-4065
USA
E-MAIL: jtaylor@stat.stanford.edu